%% file: Main.tex
\documentclass[reqno]{amsart}

\usepackage{amsmath}
\usepackage{amssymb}
\usepackage{amsfonts}
\usepackage{amsthm}
\usepackage[foot]{amsaddr}
\usepackage{mathtools}
\mathtoolsset{%
}

\usepackage[utf8]{inputenc}
\usepackage[T1]{fontenc}

\usepackage[sf,mono=false]{libertine}
\usepackage{dsfont}

\usepackage[%
cal=cm,
]
{mathalfa}

\usepackage[dvipsnames,svgnames]{xcolor}
\colorlet{MyBlue}{DodgerBlue!60!Black}
\colorlet{MyGreen}{DarkGreen!85!Black}

\usepackage[font=small,labelfont=bf]{caption}
\usepackage{subfigure}
\usepackage{tikz}
\usetikzlibrary{arrows,calc,patterns}

\usepackage{acronym}
\usepackage{latexsym}
\usepackage{nicefrac}
\usepackage{paralist}
\usepackage{wasysym}
\usepackage{xspace}

\usepackage[authoryear,compress]{natbib}

\newcommand{\citef}{\cite}
\newcommand{\citefp}{\citep}
\newcommand{\citepos}[1]{\citeauthor{#1}'s~\textpar{\citeyear{#1}}}

\usepackage{hyperref}
\hypersetup{
colorlinks=true,
linktocpage=true,
pdfstartview=FitH,
breaklinks=true,
pdfpagemode=UseNone,
pageanchor=true,
pdfpagemode=UseOutlines,
plainpages=false,
bookmarksnumbered,
bookmarksopen=false,
bookmarksopenlevel=1,
hypertexnames=true,
pdfhighlight=/O,
urlcolor=MyBlue!60!black,linkcolor=MyBlue!70!black,citecolor=DarkGreen!70!black, 
pdftitle={},
pdfauthor={},
pdfsubject={},
pdfkeywords={},
pdfcreator={pdfLaTeX},
pdfproducer={LaTeX with hyperref}
}

\numberwithin{equation}{section}  
\usepackage[sort&compress,capitalize,nameinlink]{cleveref}

\crefrangeformat{equation}{\upshape(#3#1#4)\textendash(#5#2#6)}

\newcommand{\dd}{\:d}

\newcommand{\eps}{\varepsilon}
\newcommand{\from}{\colon}

\newcommand{\pd}{\partial}

\newcommand{\ml}{\prec}

\newcommand{\R}{\mathbb{R}}

\newcommand{\N}{\mathbb{N}}

\DeclareMathOperator*{\argmax}{arg\,max}
\DeclareMathOperator*{\argmin}{arg\,min}

\DeclareMathOperator*{\union}{\bigcup}

\DeclareMathOperator{\bd}{bd}
\DeclareMathOperator{\bigoh}{\mathcal{O}}

\DeclareMathOperator{\cl}{cl}

\DeclareMathOperator{\dist}{dist}
\DeclareMathOperator{\dom}{dom}

\DeclareMathOperator{\ex}{\mathbb{E}}
\DeclareMathOperator{\grad}{\nabla}

\DeclareMathOperator{\im}{im}

\DeclareMathOperator{\intr}{int}
\DeclareMathOperator{\one}{\mathds{1}}

\DeclareMathOperator{\prob}{\mathbb{P}}

\DeclareMathOperator{\relint}{ri}

\DeclareMathOperator{\supp}{supp}

\providecommand\given{} 

\DeclarePairedDelimiter{\braces}{\{}{\}}

\DeclarePairedDelimiter{\parens}{(}{)}

\DeclarePairedDelimiter{\abs}{\lvert}{\rvert}
\DeclarePairedDelimiter{\norm}{\lVert}{\rVert}
\DeclarePairedDelimiter{\dnorm}{\lVert}{\rVert_{\ast}}

\DeclarePairedDelimiterX{\braket}[2]{\langle}{\rangle}{#1,#2}

\DeclarePairedDelimiterX{\inner}[2]{\langle}{\rangle}{#1,#2}
\DeclarePairedDelimiterX{\setdef}[2]{\{}{\}}{#1:#2}

\DeclarePairedDelimiterXPP{\probof}[1]{\prob}{(}{)}{}{%
\renewcommand\given{\nonscript\,\delimsize\vert\nonscript\,\mathopen{}}
#1}

\DeclarePairedDelimiterXPP{\exof}[1]{\ex}{[}{]}{}{%
\renewcommand\given{\nonscript\,\delimsize\vert\nonscript\,\mathopen{}}
#1}

\newcommand{\txs}{\textstyle}
\newcommand{\textpar}[1]{\textup(#1\textup)}

\newcommand{\as}{\textpar{\textrm{a.s.}}\xspace}

\usepackage[textwidth=30mm]{todonotes}
\newcommand{\debug}[1]{#1}

\newcommand{\start}{\debug 1}
\newcommand{\run}{\debug n}
\newcommand{\iRun}{\debug k}

\usepackage{algorithm}
\usepackage{algorithmic}

\theoremstyle{plain}
\newtheorem{theorem}{Theorem}
\newtheorem{corollary}[theorem]{Corollary}
\newtheorem*{corollary*}{Corollary}
\newtheorem{lemma}[theorem]{Lemma}
\newtheorem{proposition}[theorem]{Proposition}

\theoremstyle{definition}
\newtheorem{definition}[theorem]{Definition}
\newtheorem*{definition*}{Definition}

\newenvironment{Proof}[1][Proof]{\begin{proof}[#1]}{\end{proof}}

\theoremstyle{remark}
\newtheorem{remark}{Remark}
\newtheorem*{remark*}{Remark}
\newtheorem*{notation*}{Notational remark}
\newtheorem{example}{Example}

\numberwithin{theorem}{section}
\numberwithin{remark}{section}
\numberwithin{example}{section}

\newcommand{\vecspace}{\mathcal{\debug V}}
\newcommand{\dspace}{\vecspace^{\ast}}

\newcommand{\bvec}{e}

\newcommand{\cvx}{\mathcal{C}}
\newcommand{\feas}{\mathcal{\debug X}}
\newcommand{\intfeas}{\feas^{\circ}}

\newcommand{\dual}{\mathcal{\debug Y}}

\newcommand{\base}{\debug p}
\newcommand{\notbase}{\debug x}
\newcommand{\baseset}{\mathcal{C}}

\newcommand{\sol}{{\debug x}^{\ast}}

\newcommand{\tcone}{\mathrm{TC}}

\newcommand{\pcone}{\mathrm{PC}}

\newcommand{\breg}{\debug D}
\newcommand{\mirror}{\debug Q}
\newcommand{\depth}{\debug \Omega}
\newcommand{\fench}{\debug F}

\newcommand{\choice}{\mirror}
\DeclareMathOperator{\Eucl}{\Pi}
\DeclareMathOperator{\logit}{\Lambda}

\newcommand{\game}{\mathcal{\debug G}}
\newcommand{\fingame}{\debug\Gamma}

\newcommand{\play}{\debug i}
\newcommand{\playalt}{\debug j}
\newcommand{\nPlayers}{\debug N}
\newcommand{\players}{\mathcal{\nPlayers}}

\newcommand{\playOne}{\debug A}
\newcommand{\playTwo}{\debug B}

\newcommand{\pure}{\debug\alpha}
\newcommand{\purealt}{\debug\beta}
\newcommand{\nPures}{\debug A}
\newcommand{\pures}{\mathcal{\debug\nPures}}
\newcommand{\peq}{\pure^{\ast}}

\newcommand{\act}{\debug X}
\newcommand{\actions}{\feas}
\newcommand{\acts}{\actions}
\newcommand{\intacts}{\intfeas}

\newcommand{\pay}{\debug u}
\newcommand{\payv}{\debug v}

\newcommand{\cost}{c}
\newcommand{\pot}{f}

\newcommand{\est}{\hat\payv}
\newcommand{\score}{\debug Y}

\newcommand{\gamefull}{\game(\players,(\acts_{\play})_{\play\in\players},(\pay_{\play})_{\play\in\players})}
\newcommand{\fingamefull}{\fingame(\players,(\pures_{\play})_{\play\in\players},(\pay_{\play})_{\play\in\players})}

\newcommand{\eq}{\sol}
\newcommand{\eqset}{\feas^{\ast}}
\newcommand{\eqnhd}{U^{\ast}}
\newcommand{\payveq}{\payv^{\ast}}
\newcommand{\deq}{y^{\ast}}

\newcommand{\dev}{\debug q}
\newcommand{\val}{\pay^{\ast}}
\newcommand{\gap}{\debug\epsilon}
\newcommand{\length}{\ell}

\newcommand{\reg}{\debug R}

\newcommand{\filter}{\mathcal{F}}

\newcommand{\simplex}{\Delta}
\newcommand{\intsimplex}{\simplex^{\!\circ}}

\newcommand{\step}{\debug\gamma}

\newcommand{\resource}{r}
\newcommand{\nResources}{R}
\newcommand{\resources}{\mathcal{\nResources}}
\newcommand{\load}{w}

\newcommand{\route}{\pure}
\newcommand{\routes}{\pures}

\newcommand{\hessmat}{H^{\game}}
\newcommand{\vbound}{\debug V_{\ast}}
\newcommand{\noise}{\debug\xi}
\newcommand{\noisedev}{\debug\sigma}
\newcommand{\noisevar}{\noisedev^{2}}
\newcommand{\snoise}{\psi}
\newcommand{\semiflow}{\debug\Phi}

\newcommand{\strong}{\debug K}
\renewcommand{\sharp}{\debug c}

\newcommand{\ie}{i.e.,\xspace}
\newcommand{\eg}{e.g.,\xspace}

\begin{document}


\title
[Learning in games with continuous action sets]
{Learning in games with continuous action sets\\and unknown payoff functions}

\author
[P.~Mertikopoulos]
{Panayotis Mertikopoulos$^{1}$}
\address{$^{1}$
Univ. Grenoble Alpes, CNRS, Inria, LIG, F-38000, Grenoble, France.}
\email{\href{mailto:panayotis.mertikopoulos@imag.fr}{panayotis.mertikopoulos@imag.fr}}
\author
[Z.~Zhou]
{Zhengyuan Zhou$^{2}$}
\address{$^{2}$
Stanford University, Dept. of Electrical Engineering, Stanford, CA, 94305.}
\email{\href{mailto:zyzhou@stanford.edu}{zyzhou@stanford.edu}}

\input{Thanks}

\subjclass[2010]{Primary 91A26, 90C15; secondary 90C33, 68Q32.}
\keywords{%
Continuous games;
dual averaging;
variational stability;
Fenchel coupling;
Nash equilibrium.}

\newcommand{\acdef}[1]{\textit{\acl{#1}} \textup{(\acs{#1})}\acused{#1}}
\newcommand{\acdefp}[1]{\emph{\aclp{#1}} \textup(\acsp{#1}\textup)\acused{#1}}
\newcommand{\acli}[1]{\textit{\acl{#1}}}

\newacro{CCE}{coarse correlated equilibrium}
\newacroplural{CCE}[CCE]{coarse correlated equilibria}
\newacro{CE}{correlated equilibrium}
\newacroplural{CE}[CE]{correlated equilibria}
\newacro{OMD}{online mirror descent}
\newacro{OGD}{online gradient descent}
\newacro{MDS}{martingale difference sequence}
\newacro{LLN}{law of large numbers}
\newacro{EW}{exponential weights}
\newacro{APT}{asymptotic pseudotrajectory}
\newacroplural{APT}{asymptotic pseudotrajectories}
\newacro{KKT}{Karush\textendash Kuhn\textendash Tucker}
\newacro{MSE}{mean squared error}
\newacro{DSC}{diagonal strict concavity}
\newacro{NE}{Nash equilibrium}
\newacroplural{NE}[NE]{Nash equilibria}
\newacro{SD}[s/d]{source/destination}
\newacro{ESS}{evolutionarily stable state}
\newacro{MD}{mirror descent}
\newacro{MA}{mirror ascent}
\newacro{MMA}{multi-player mirror ascent}
\newacro{DA}{dual averaging}
\newacro{MDA}{multi-agent dual averaging}
\newacro{SA}{stochastic approximation}
\newacro{VI}{va\-ri\-a\-tio\-nal inequality}
\newacroplural{VI}{va\-ri\-a\-tio\-nal inequalities}
\newacro{VS}{variational stability}
\newacro{ES}{evolutionary stability}
\newacro{GNEP}{generalized Nash equilibrium problem}
\newacro{iid}[i.i.d.]{independent and identically distributed}

\begin{abstract}
\input{Abstract}
\end{abstract}
\acresetall

\maketitle



\section{Introduction}
\label{sec:introduction}
\input{Introduction}

\section{Continuous games and variational stability}
\label{sec:prelims}
\input{Prelims}

\section{Learning via dual averaging}
\label{sec:learning}
\input{Learning}

\section{Convergence analysis}
\label{sec:analysis}
\input{Analysis}

\section{Learning in finite games}
\label{sec:finite}
\input{Finite}

\section{Speed of convergence}
\label{sec:rates}
\input{Rates}

\section{Discussion}
\label{sec:discussion}
\input{Discussion}

\appendix


\section{Auxiliary results}
\label{app:aux}
\input{App-Auxiliary}

\bibliographystyle{apalike}
\bibliography{IEEEabrv,Bibliography}

\end{document}

%% file: Thanks.tex
%
%
\thanks{%
The authors are indebted to the associate editor and two anonymous referees for their detailed suggestions and remarks.
The paper has also benefited greatly from thoughtful comments by Jérôme Bolte, Nicolas Gast, Jérôme Malick, Mathias Staudigl, and the audience of the Paris Optimization Seminar.}
\thanks{%
P.~Mertikopoulos was partially supported by
the French National Research Agency (ANR) project ORACLESS (ANR\textendash GAGA\textendash13\textendash JS01\textendash 0004\textendash 01)
and the Huawei Innovation Research Program ULTRON}

%% file: Abstract.tex
%
%
This paper examines the convergence of no-regret learning in games with continuous action sets.
For concreteness, we focus on learning via ``dual averaging'', a widely used class of no-regret learning schemes where players take small steps along their individual payoff gradients and then ``mirror'' the output back to their action sets.
In terms of feedback, we assume that players can only estimate their payoff gradients up to a zero-mean error with bounded variance.
To study the convergence of the induced sequence of play, we introduce the notion of \emph{\acl{VS}}, and we show that stable equilibria are locally attracting with high probability whereas globally stable equilibria are globally attracting with probability $1$.
We also discuss some applications to mixed-strategy learning in finite games, and we provide explicit estimates of the method's convergence speed.

%% file: Introduction.tex

The prototypical setting of online optimization can be summarized as follows:
at every stage $\run = 1, 2,\dotsc$, of a repeated decision process, an agent selects an action $\act_{\run}$ from some set $\feas$ (assumed here to be convex and compact), and obtains a reward $\pay_{\run}(\act_{\run})$ determined by an \emph{a priori} unknown payoff function $\pay_{\run}\from\feas\to\R$.
Subsequently, the agent receives some problem-specific feedback (for instance, an estimate of the gradient of $\pay_{\run}$ at $\act_{\run}$), and selects a new action with the goal of maximizing the obtained reward.
Aggregating over the stages of the process, this is usually quantified by asking that the agent's \emph{regret} $\reg_{\run} \equiv \max_{x\in\feas} \sum_{\iRun=1}^{\run} \left[ \pay_{\iRun}(x) - \pay_{\iRun}(\act_{\iRun}) \right]$ grow sublinearly in $\run$, a property known as ``no regret''.

In this general setting, the most widely used class of no-regret policies is the \acdef{OMD} method of \citef{SS07} and its variants \textendash\ such as ``Following the Regularized Leader'' \citep{SSS07}, \acl{DA} \citep{Nes09,Xia10}, etc.
Specifically, if the problem's payoff functions are concave, \acl{MD} guarantees an $\bigoh(\sqrt{\run})$ regret bound which is well-known to be tight in a ``black-box'' environment (\ie without any further assumptions on $\pay_{\run}$).
Thus, owing to these guarantees, this class of first-order methods has given rise to an extensive literature in online learning and optimization;
for a survey, see \cite{SS11}, \cite{BCB12}, \cite{Haz12}, and references therein.

In this paper, we consider a multi-agent extension of the above framework where the agents' rewards are determined by their individual actions and the actions of all other agents via a fixed mechanism:
a \emph{non-cooperative game}.
Even though this mechanism may be unknown and/or opaque to the players, the additional structure it provides means that finer convergence criteria apply, chief among them being that of convergence to a \acdef{NE}.
We are thus led to the following fundamental question:
\emph{if all players of a repeated game employ a no-regret updating policy, do their actions converge to a \acl{NE} of the underlying game?}

\subsection*{Summary of contributions}

In general, the answer to this question is a resounding ``no''.
Even in simple, finite games, no-regret learning may cycle \citep{MPP18} and its limit set may contain highly non-rationalizable strategies that assign positive weight \emph{only} to strictly dominated strategies \citep{VZ13}.
As such, our aim in this paper is twofold:
\begin{enumerate}
[\indent\itshape i\hspace*{1pt}\upshape)]
\item
to provide sufficient conditions under which no-regret learning converges to equilibrium;
and
\item
to assess the speed and robustness of this convergence in the presence of uncertainty, feedback noise, and other learning impediments.
\end{enumerate}

Our contributions along these lines are as follows:
First, in \cref{sec:prelims}, we introduce an equilibrium stability notion which we call \acdef{VS}, and which is formally similar to (and inspired by) the seminal notion of \acli{ES} in population games \citep{MSP73}.%
\footnote{Heuristically, \acl{VS} is to games with a finite number of players and a continuum of actions what \acl{ES} is to games with a continuum of players and a finite action space.
Our choice of terminology reflects precisely this analogy.}
This stability notion extends the standard notion of operator monotonicity, so it applies in particular to all \emph{monotone games} (that is, concave games that satisfy \citepos{Ros65} \acl{DSC} condition).
In fact, going beyond concave games, \acl{VS} allows us to treat convergence questions in general games with continuous action spaces \emph{without} having to restrict ourselves to a specific subclass (such as potential or common interest games).

Our second contribution is a detailed analysis of the long-run behavior of no-regret learning under \acl{VS}.
Regarding the information available to the players, our only assumption is that they have access to unbiased, bounded-variance estimates of their individual payoff gradients at each step;
beyond this, we assume no prior knowledge of their payoff functions and/or the game.
Despite this lack of information, \acl{VS} guarantees that
\begin{inparaenum}
[(\itshape i\hspace*{1pt}\upshape)]
\item
the induced sequence of play converges globally to globally stable equilibria with probability~$1$ (\cref{thm:global-imperfect});
and
\item
it converges locally to locally stable equilibria with high probability (\cref{thm:local-imperfect}).
\end{inparaenum}
As a corollary, if the game admits a \textpar{pseudo-}concave potential or if it is monotone, the players' actions converge to \acl{NE} no matter the level of uncertainty affecting the players' feedback.
In \cref{sec:finite}, we further extend these results to learning with imperfect feedback in finite games.

Our third contribution concerns the method's convergence speed.
Mirroring a known result of \cite{Nes09} for \aclp{VI}, we show that the gap from a stable state decays ergodically as $\bigoh(1/\sqrt{\run})$ if the method's step-size is chosen appropriately.
Dually to this, we also show that the algorithm's expected running length until players reach an $\eps$-neighborhood of a stable state is $\bigoh(1/\eps^{2})$.
Finally, if the stage game admits a \emph{sharp} equilibrium (a straightforward extension of the notion of strict equilibrium in finite games), we show that, with probability $1$, the process reaches an equilibrium in a \emph{finite} number of steps.

\color{black}
Our analysis relies on tools and techniques from stochastic approximation, martingale limit theory and convex analysis.
In particular, with regard to the latter, we make heavy use of a ``primal-dual divergence'' measure between action and gradient variables, which we call the \emph{Fenchel coupling}.
This coupling is a hybridization of the Bregman divergence which provides a potent tool for proving convergence thanks to its Lyapunov properties.

\subsection*{Related work}
Originally, mirror descent was introduced by \citef{NY83} for solving offline convex programs.
The \acdef{DA} variant that we consider here was pioneered by \cite{Nes09} and proceeds as follows:%
\footnote{In the online learning literature, \acl{DA} is sometimes called lazy \acl{MD} and can be seen as a linearized ``Follow the Regularized Leader'' (FTRL) scheme \textendash\ for more details, we refer the reader to \cite{BecTeb03}, \cite{Xia10}, and \cite{SS11}.}
at each stage, the method takes a gradient step in a dual space (where gradients live);
the result is then mapped (or ``mirrored'') back to the problem's feasible region, a new gradient is generated, and the process repeats.
The ``mirroring'' step above is itself determined by a strongly convex regularizer (or ``distance generating'') function:
the squared Euclidean norm gives rise to \citepos{Zin03} \acl{OGD} algorithm,
while the (negative) Gibbs entropy on the simplex induces the well-known \ac{EW} algorithm \citep{Vov90,AHK12}.

\cite{Nes09} and \cite{NJLS09} provide several convergence results for \acl{DA} in (stochastic) convex programs and saddle-point problems, while \cite{Xia10} provides a thorough regret analysis for online optimization problems.
In addition to treating the interactions of several competing agents at once, the fundamental difference of our paper with these works is that the convergence analysis in the latter is ``ergodic'', \ie it concerns the time-averaged sequence $\bar\act_{\run} = \sum_{\iRun=1}^{\run} \step_{\iRun} \act_{\iRun} / \sum_{\iRun=1}^{\run} \step_{\iRun}$, and \emph{not} the actual sequence of actions $\act_{\run}$ employed by the players.

In online optimization, this averaging comes up naturally because the focus is on the players' regret.
In the offline case, the points where an oracle is called during the execution of an algorithm do not carry any particular importance, so averaging provides a convenient way of obtaining convergence.
However, in a game-theoretic setting, the figure of merit is the \emph{actual sequence of play}, which determines the players' payoffs at each stage.
The behavior of $\act_{\run}$ may differ drastically from that of $\bar \act_{\run}$, so our treatment requires a completely different set of tools and techniques (especially in the stochastic regime).

Much of our analysis boils down to solving in an online way a (stochastic) \ac{VI} characterizing the game's \aclp{NE}.
\cite{Nes07} and \cite{JNT11} provide efficient offline methods to do this, relying on an ``extra-gradient'' step to boost the convergence rate of the ergodic sequence $\bar \act_{\run}$.
In our limited-feedback setting, we do not assume that players can make an extra oracle call to actions that were not actually employed, so the extrapolation results of \cite{Nes07} and \cite{JNT11} do not apply.
The single-call results of \cite{Nes09} are closer in spirit to our paper but, again, they focus exclusively on monotone \aclp{VI} and the ergodic sequence $\bar \act_{\run}$ \textendash\ not the actual sequence of play $\act_{\run}$.
All the same, for completeness, we make the link with ergodic convergence in \cref{thm:zerosum,thm:gap-ergodic}.

When applied to mixed-strategy learning in finite games, the class of algorithms studied here has very close ties to the family of perturbed best response maps that arise in models of fictitious play and reinforcement learning \citep{HS02,LC05,CGM15}.
Along these lines, \citef{MS16} recently showed that a continuous-time version of the dynamics studied in this paper eliminates dominated strategies and converges to strict equilibria from all nearby initial conditions.
Our analysis in \cref{sec:finite} extends these results to a discrete-time, stochastic setting.

In games with continuous action sets, \citef{PL12} and \citef{PML17} examined a mixed-strategy actor-critic algorithm which converges to a probability distribution that assigns most weight to equilibrium states.
At the pure strategy level, several authors have considered \ac{VI}-based and Gauss\textendash Seidel methods for solving \acp{GNEP};
for a survey, see \citef{FK07} and \citef{SFPP10}.
The intersection of these works with the current paper is when the game satisfies a global monotonicity condition similar to the \acl{DSC} condition of \citef{Ros65}.
However, the literature on \acp{GNEP} does not consider the implications for the players' regret, the impact of uncertainty and/or local convergence/stability issues, so there is no overlap with our results.

Finally, during the final preparation stages of this paper (a few days before the actual submission), we were made aware of a preprint by \citef{BBF16} examining the convergence of pure-strategy learning in strictly concave games with one-di\-men\-sio\-nal action sets.
A key feature of the analysis of \citef{BBF16} is that players only observe their realized, in-game payoffs, and they choose actions based on their payoffs' variation from the previous period. 
The resulting mean dynamics boil down to an instantiation of \acl{DA} induced by the entropic regularization penalty $h(x) = x\log x$ (cf. \cref{sec:learning}), suggesting several interesting links with the current work.

\subsection*{Notation}
\label{sec:notation}

Given a finite-dimensional vector space $\vecspace$ with norm $\norm{\cdot}$,
we write
$\dspace$ for its dual,
$\braket{y}{x}$ for the pairing between $y\in\dspace$ and $x\in\vecspace$,
and $\dnorm{y} \equiv \sup\setdef{\braket{y}{x}}{\norm{x} \leq 1}$ for the dual norm of $y$ in $\dspace$.
If $\cvx\subseteq\vecspace$ is convex, we also write $\cvx^{\circ} \equiv \relint(\cvx)$ for the relative interior of $\cvx$, $\norm{\cvx} = \sup\setdef{\norm{x' - x}}{x,x'\in\cvx}$ for its diameter, and $\dist(\cvx,x) = \inf_{x'\in\cvx} \norm{x' - x}$ for the distance between $x\in\vecspace$ and $\cvx$.

For a given $x\in\cvx$, the \emph{tangent cone} $\tcone_{\cvx}(x)$ is defined as the closure of the set of all rays emanating from $x$ and intersecting $\cvx$ in at least one other point;
dually, the \emph{polar cone} $\pcone_{\cvx}(x)$ to $\cvx$ at $x$ is defined as $\pcone_{\cvx}(x) = \setdef{y\in\dspace}{\braket{y}{z} \leq 0\; \text{for all}\; z\in\tcone_{\cvx}(x)}$.
For concision, when $\cvx$ is clear from the context, we will drop it altogether and write $\tcone(x)$ and $\pcone(x)$ instead.

%% file: Prelims.tex

\subsection{Basic definitions and examples}
\label{sec:games}

Throughout this paper, we focus on games played by a finite set of \emph{players} $\play\in\players = \{1,\dotsc,\nPlayers\}$.
During play, each player selects an \emph{action} $x_{\play}$ from a compact convex subset $\feas_{\play}$ of a finite-dimensional normed space $\vecspace_{\play}$,
and
their reward is determined by the profile $x=(x_{1},\dotsc,x_{\nPlayers})$ of all players' actions \textendash\ often denoted as $x \equiv (x_{\play};x_{-\play})$ when we seek to highlight the action $x_{\play}$ of player $\play$ against the ensemble of actions $x_{-\play} = (x_{\playalt})_{\playalt\neq\play}$ of all other players.

In more detail, writing $\feas \equiv \prod_{\play} \feas_{i}$ for the game's \emph{action space}, each player's \emph{payoff} is determined by an associated \emph{payoff function} $\pay_{\play}\from\feas\to\R$.
In terms of regularity, we assume that $\pay_{\play}$ is continuously differentiable in $x_{\play}$, and we write
\begin{equation}
\label{eq:grad}
\payv_{\play}(x)
	\equiv \nabla_{x_{\play}} \pay_{\play}(x_{\play};x_{-\play})
\end{equation}
for the \emph{individual gradient} of $\pay_{\play}$ at $x$;
we also assume that $\pay_{\play}$ and $\payv_{\play}$ are both continuous in $x$.%
\footnote{\label{foot:dual}%
In the above, we tacitly assume that $\pay_{\play}$ is defined on an open neighborhood of $\feas_{\play}$.
This allows us to use ordinary derivatives, but none of our results depend on this device.
We also note that $\payv_{\play}(x)$ acts naturally on vectors $z_{\play}\in\vecspace_{\play}$ via the mapping $z_{\play} \mapsto \braket{\payv_{\play}(x)}{z_{\play}} \equiv \pay_{\play}'(x;z_{i}) = d/d\tau\vert_{\tau=0}\, \pay_{\play}(x_{\play}+\tau z_{\play};x_{-\play})$;
in view of this, $\payv_{\play}(x)$ is treated as an element of $\dspace_{\play}$, the dual of $\vecspace_{\play}$.}
Putting all this together, a \emph{continuous game} is a tuple $\game \equiv \gamefull$ with players, actions and payoffs defined as above.

\smallskip

As a special case, we will sometimes consider payoff functions that are \emph{individually \textpar{pseudo-}concave} in the sense that
\begin{equation}
\label{eq:concavity}
\txs
\text{\itshape
$\pay_{\play}(x_{\play};x_{-\play})$ is \textpar{pseudo-}concave in $x_{\play}$
	for all $x_{-\play}\in\prod_{\playalt\neq\play}\feas_{\playalt}$,
	$\play\in\players$.
}
\end{equation}
When this is the case, we say that the game itself is \textpar{pseudo-}concave.
Below, we briefly discuss some well-known examples of such games:

\begin{example}[Mixed extensions of finite games]
\label{ex:finite}
In a \emph{finite game} $\fingame\equiv(\players,\pures,\pay)$, each player $\play\in\players$ chooses an action $\pure_{\play}$ from a finite set $\pures_{\play}$ of ``pure strategies'' and no assumptions are made on the players' payoff functions $\pay_{\play}\from\pures\equiv\prod_{\playalt}\pures_{\playalt}\to\R$.
Players can ``mix'' these choices by playing \emph{mixed strategies},
\ie probability distributions $x_{\play}$ drawn from the simplex $\feas_{\play} \equiv \simplex(\pures_{\play})$.
In this case (and in a slight abuse of notation), the expected payoff to player $\play$ in the mixed profile $x = (x_{1},\dotsc,x_{\nPlayers})$ can be written as
\begin{equation}
\label{eq:pay-mixed}
\pay_{\play}(x)
	= \sum_{\pure_{1}\in\pures_{1}} \dotsm\;\sum_{\crampedclap{\pure_{\nPlayers}\in\pures_{\nPlayers}}} \;
	\pay_{\play}(\pure_{1},\dotsc,\pure_{\nPlayers})\;
	x_{1,\pure_{1}} \dotsm\, x_{\nPlayers,\pure_{\nPlayers}},
\end{equation}
so the players' individual gradients are simply their payoff vectors:
\begin{equation}
\label{eq:payv-finite}
\payv_{\play}(x)
	= \nabla_{x_{\play}} \pay_{\play}(x)
	= (\pay_{\play}(\pure_{\play};x_{-\play}))_{\pure_{\play}\in\pures_{\play}}.
\end{equation}
The resulting continuous game is called the \emph{mixed extension} of $\fingame$.
Since $\feas_{\play} = \simplex(\pures_{\play})$ is convex and $\pay_{\play}$ is linear in $x_{\play}$, $\game$ is itself concave in the sense of \eqref{eq:concavity}.
\end{example}

\begin{example}[Cournot competition]
\label{ex:Cournot}
Consider the following Cournot o\-li\-go\-po\-ly model:
There is a finite set $\players = \{1,\dotsc,\nPlayers\}$ of \emph{firms}, each supplying the market with a quantity $x_{\play} \in [0,C_{\play}]$ of the same good (or service) up to the firm's production capacity $C_{\play}$.
This good is then priced as a decreasing function $P(x)$ of each firm's production;
for concreteness, we focus on the linear model $P(x) = a - \sum_{\play} b_{\play} x_{\play}$
where $a$ is a positive constant and the coefficients $b_{\play}>0$ reflect the price-setting power of each firm.

In this model, the utility of firm $\play$ is given by
\begin{equation}
\label{eq:pay-Cournot}
\pay_{\play}(x)
	= x_{\play} P(x) - c_{\play} x_{\play},
\end{equation}
where $c_{\play}$ represents the marginal production cost of firm $\play$.
Letting $\feas_{\play} = [0,C_{\play}]$, the resulting game is easily seen to be concave in the sense of \eqref{eq:concavity}.
\end{example}

\begin{example}[Congestion games]
\label{ex:congestion}
Congestion games are game-theoretic models that arise in the study of traffic networks (such as the Internet).
To define them, fix a set of players $\players$ that share a set of \emph{resources} $\resource\in\resources$,
each associated with a nondecreasing convex \emph{cost function} $\cost_{\resource}\from\R_{+}\to\R$ (for instance, links in a data network and their corresponding delay functions).
Each player $\play\in\players$ has a certain \emph{resource load} $\rho_{\play} > 0$ which is split over a collection $\routes_{\play}\subseteq 2^{\resources}$ of resource subsets $\route_{\play}$ of $\resources$ \textendash\ \eg sets of links that form paths in the network.
Then, the action space of player $\play\in\players$ is the scaled simplex $\feas_{\play} = \rho_{\play} \simplex(\routes_{\play}) = \setdef{x_{\play} \in \R_{+}^{\routes_{\play}}}{\sum_{\route_{\play}\in\routes_{\play}} x_{\play\route_{\play}} = \rho_{\play}}$ of \emph{load distributions} over $\routes_{\play}$.

Given a load profile $x = (x_{1},\dotsc,x_{\nPlayers})$, costs are determined based on the utilization of each resource as follows:
First, the \emph{demand} $\load_{\resource}$ of the $\resource$-th resource is defined as the total load $\load_{\resource} = \sum_{\play\in\players} \sum_{\route_{\play}\ni\resource} x_{\play\route_{\play}}$ on said resource.
This demand incurs a \emph{cost} $\cost_{\resource}(\load_{\resource})$ per unit of load to each player utilizing resource $\resource$, where $\cost_{\resource}\from\R_{+}\to\R$ is a nondecreasing convex function.
Accordingly, the total cost to player $\play\in\players$ is
\begin{equation}
\label{eq:pay-congestion}
\cost_{\play}(x)
	= \sum_{\route_{\play}\in\routes_{\play}}
	x_{\play\route_{\play}} \cost_{\play\route_{\play}}(x),
\end{equation}
where $\cost_{\play\route_{\play}}(x) = \sum_{\resource\in\route_{\play}} \cost_{\resource}(\load_{\resource})$ denotes the cost incurred to player $\play$ by the utilization of $\route_{\play}\subseteq\resources$.
The resulting \emph{atomic splittable congestion game} $\game\equiv\game(\players,\feas,-\cost)$ is easily seen to be concave in the sense of \eqref{eq:concavity}.
\end{example}

\subsection{\acl{NE}}
\label{sec:Nash}

Our analysis focuses primarily on \acdefp{NE}, \ie strategy profiles that discourage unilateral deviations.
Formally, $\eq\in\feas$ is a \acli{NE}
if
\begin{equation*}
\label{eq:Nash}
\tag{NE}
\pay_{\play}(\eq_{\play};\eq_{-\play})
	\geq \pay_{\play}(x_{\play};\eq_{-\play})
	\quad
	\text{for all $x_{\play}\in\feas_{\play}$, $\play\in\players$}.
\end{equation*}
Obviously, if $\eq$ is a \acl{NE}, we have the first-order condition
\begin{equation}
\label{eq:Nash-vec}
\pay_{\play}'(\eq;z_{\play})
	= \braket{\payv_{\play}(\eq)}{z_{\play}}
	\leq 0
	\quad
	\text{for all $z_{\play}\in\tcone_{\play}(\eq_{\play})$, $\play\in\players$},
\end{equation}
where $\tcone_{\play}(\eq_{\play})$ denotes the \emph{tangent cone} to $\feas_{\play}$ at $\eq_{\play}$.
Therefore, if $\eq$ is a \acl{NE}, each player's individual gradient $\payv_{\play}(\eq)$ belongs to the \emph{polar cone} $\pcone_{\play}(\eq_{\play})$ to $\feas_{\play}$ at $\eq_{\play}$ (cf.~\cref{fig:Nash});
moreover,
the converse also holds if the game is pseudo-concave.
We encode this more concisely as follows:


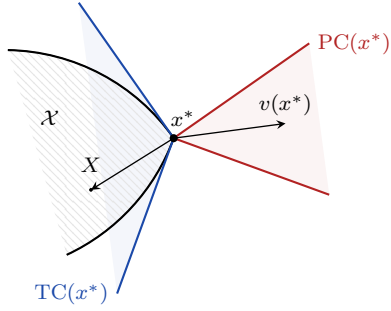
\begin{figure}
\centering
\footnotesize
\input{Figures/Nash}
\caption{Geometric characterization of \aclp{NE}.}
\label{fig:Nash}
\end{figure}


\begin{proposition}
\label{prop:Nash-var}
If $\eq\in\feas$ is a \acl{NE}, then $\payv(\eq)\in\pcone(\eq)$, \ie
\begin{equation}
\label{eq:Nash-var}
\braket{\payv(\eq)}{x - \eq}
	\leq 0
	\quad
	\text{for all $x\in\feas$}.
\end{equation}
The converse also holds if the game is \textpar{pseudo-}concave in the sense of \eqref{eq:concavity}.
\end{proposition}

\begin{remark}
In the above (and in what follows), $\payv = (\payv_{\play})_{\play\in\players}$ denotes the ensemble of the players' individual payoff gradients
and $\braket{\payv}{z} \equiv \sum_{\play\in\players} \braket{\payv_{\play}}{z_{\play}}$ stands for the pairing between $\payv$ and the vector $z = (z_{\play})_{\play\in\players} \in \prod_{\play\in\players} \vecspace_{\play}$.
For concision, we also write $\vecspace \equiv \prod_{\play} \vecspace_{\play}$ for the ambient space of $\feas\equiv\prod_{\play}\feas_{\play}$ and $\dspace$ for its dual.
\end{remark}

\begin{Proof}[Proof of \cref{prop:Nash-var}]
If $\eq$ is a \acl{NE}, \eqref{eq:Nash-var} is obtained by setting $z_{\play} = x_{\play} - \eq_{\play}$ in \eqref{eq:Nash-vec} and summing over all $\play\in\players$.
Conversely, if \eqref{eq:Nash-var} holds and the game is \textpar{pseudo-}concave, pick some $x_{\play}\in\feas_{\play}$
and let $x = (x_{\play};\eq_{-\play})$ in \eqref{eq:Nash-var}.
This gives $\braket{\payv_{\play}(\eq)}{x_{\play} - \eq_{\play}} \leq 0$ for all $x_{\play}\in\feas_{\play}$ so \eqref{eq:Nash} follows by the basic properties of (pseudo-)concave functions.
\end{Proof}

\cref{prop:Nash-var} shows that \aclp{NE} of concave games are precisely the solutions of the \acl{VI} \eqref{eq:Nash-var}, so existence follows from standard results.
Using a similar variational characterization, \cite{Ros65} proved the following sufficient condition for equilibrium uniqueness:

\begin{theorem}[\citeauthor{Ros65}, \citeyear{Ros65}]
\label{thm:Rosen}
Assume that $\game$ satisfies the payoff monotonicity condition
\begin{equation*}
\label{eq:MC}
\tag{MC}
\braket{\payv(x') - \payv(x)}{x' - x}
	\leq 0
	\quad
	\text{for all $x,x'\in\feas$},
\end{equation*}
with equality if and only if $x=x'$.
Then, $\game$ admits a unique \acl{NE}.
\end{theorem}

Games satisfying \eqref{eq:MC} are called \emph{\textpar{strictly} monotone} and they enjoy properties similar to those of (strictly) convex functions.%
\footnote{\cite{Ros65} originally referred to \eqref{eq:MC} as \acl{DSC};
\cite{HS09} use the term ``stable'' for population games that satisfy a formal analogue of \eqref{eq:MC}, while \cite{San15} and \cite{SW16} call such games ``contractive'' and ``dissipative'' respectively.
In all cases, the adverb ``strictly'' refers to the ``only if'' requirement in \eqref{eq:MC}.}
In particular, letting $x_{-\play}' = x_{-\play}$, \eqref{eq:MC} gives
\begin{equation}
\braket{\payv_{\play}(x_{\play}';x_{-\play}) - \payv_{\play}(x_{\play};x_{-\play})}{x_{\play}' - x_{\play}}
	\leq 0
	\quad
	\text{for all $x_{\play},x_{\play}'\in\feas_{\play}$, $x_{-\play}\in\feas_{-\play}$},
\end{equation}
implying in turn that $\pay_{\play}(x)$ is (strictly) concave in $x_{\play}$ for all $\play$.
Therefore,
\emph{any game satisfying \eqref{eq:MC} is also concave.}

\subsection{Variational stability}
\label{sec:stability}

Combining \cref{prop:Nash-var} and \eqref{eq:MC}, it follows that the (necessarily unique) \acl{NE} of a monotone game satisfies the inequality
\begin{equation}
\label{eq:VS-0}
\braket{\payv(x)}{x - \eq}
	\leq \braket{\payv(\eq)}{x - \eq}
	\leq 0
	\quad
	\text{for all $x\in\feas$}.
\end{equation}
In other words, if $\eq$ is a \acl{NE} of a monotone game, the players' individual payoff gradients ``point towards'' $\eq$ in the sense that $\payv(x)$ forms an acute angle with $\eq-x$.
Motivated by this, we introduce below the following relaxation of the monotonicity condition \eqref{eq:MC}:

\begin{definition}
\label{def:VS-point}
We say that $\eq\in\feas$ is \emph{variationally stable} \textpar{or simply \emph{stable}} if there exists a neighborhood $U$ of $\eq$ such that
\begin{equation*}
\notag
\braket{\payv(x)}{x - \eq} 
	\leq 0
	\quad
	\text{for all $x\in U$},
\end{equation*}
with equality if and only if $x=\eq$.
In particular, if $U$ can be taken to be all of $\feas$, we say that $\eq$ is \emph{globally variationally stable} \textpar{or \emph{globally stable} for short}.
\end{definition}

\begin{remark}
The terminology ``variational stability'' alludes to the seminal notion of \acli{ES} introduced by \cite{MSP73} for population games (\ie games with a continuum of players and a common, finite set of actions $\pures$).
Specifically, if $\payv(x) = (\payv_{\pure}(x))_{\pure\in\pures}$ denotes the payoff field of such a game (with $x\in\simplex(\pures)$ denoting the state of the population), \cref{def:VS} boils down to the variational characterization of \aclp{ESS} due to \cite{HSS79}.
As we show in the next sections, \acl{VS} plays the same role for learning in games with continuous action spaces as \acl{ES} plays for evolution in games with a continuum of players.
\end{remark}

By \eqref{eq:VS-0}, a first example of \acl{VS} is provided by the class of monotone games:

\begin{corollary}
\label{cor:stable}
If $\game$ satisfies \eqref{eq:MC}, its \textpar{unique} \acl{NE} is globally stable.
\end{corollary}

The converse to \cref{cor:stable} does not hold, even partially.
For instance, consider the single-player game with payoffs given by the function
\begin{equation}
\label{eq:pay-example}
\pay(x)
	= 1 - \sum_{\ell=1}^{d} \sqrt{1+x_{\ell}},
	\quad
	x\in[0,1]^{d}.
\end{equation}
In this simple example, the origin is the unique maximizer (and hence unique \acl{NE}) of $\pay$.
Moreover, we trivially have $\braket{\payv(x)}{x} = -2\sum_{\ell=1}^{d} x_{\ell}/\sqrt{1+x_{\ell}} \leq 0$ with equality if and only if $x=0$, so the origin satisfies the global version of \eqref{eq:VS};
however, $\pay$ is not even pseudo-concave if $d\geq2$, so the game cannot be monotone.
In words, \eqref{eq:MC} is a sufficient condition for the existence of a (globally) stable state, but not a necessary one.

Nonetheless, even in this (non-monotone) example, variational stability characterizes the game's unique \acl{NE}.
We make this link precise below:



\begin{table}[tbp]
\centering
\renewcommand{\arraystretch}{1.3}
\small
\input{Tables/Stability}
\vspace{2ex}
\caption{Monotonicity, stability, and \acl{NE}:
the existence of a concave potential implies monotonicity;
monotonicity implies the existence of a globally stable point;
and
globally stable points are equilibria.
}
\label{tab:stability}
\vspace{-3ex}
\end{table}


\begin{proposition}
\label{prop:Nash-stable-point}
Suppose that $\eq\in\feas$ is variationally stable.
Then:
\begin{enumerate}
[\indent a\upshape)]
\item
If $\game$ is \textpar{pseudo-}concave, $\eq$ is an isolated \acl{NE} of $\game$.
\item
If $\eq$ is globally stable, it is the game's unique \acl{NE}.
\end{enumerate}
\end{proposition}

\cref{prop:Nash-stable-point} indicates that variationally stable states are isolated (for the proof, see that of \cref{prop:Nash-stable} below).
However, this also means that \aclp{NE} of games that admit a concave \textendash\ but not \emph{strictly} concave \textendash\ potential may fail to be stable.
To account for such cases, we will also consider the following setwise version of \acl{VS}:

\begin{definition}
\label{def:VS}
Let $\eqset\subseteq\feas$ be closed and nonempty.
We say that $\eqset$ is
\emph{variationally stable} \textpar{or simply \emph{stable}} if there exists a neighborhood $U$ of $\eqset$ such that
\begin{equation}
\label{eq:VS}
\tag{VS}
\braket{\payv(x)}{x - \eq} 
	\leq 0
	\quad
	\text{for all $x\in U$, $\eq\in\eqset$},
\end{equation}
with equality for a given $\eq\in\eqset$ if and only if $x\in\eqset$.
In particular, if $U$ can be taken to be all of $\feas$, we say that $\eqset$ is \emph{globally variationally stable} \textpar{or \emph{globally stable} for short}.
\end{definition}

Obviously, \cref{def:VS} subsumes \cref{def:VS-point}:
if $\eq\in\feas$ is stable in the pointwise sense of \cref{def:VS-point}, then it is also stable when viewed as a singleton set.
In fact, when this is the case, it is also easy to see that $\eq$ cannot belong to some larger variationally stable set,%
\footnote{In that case \eqref{eq:VS} would give $\braket{\payv(x')}{x' - \eq} = 0$ for some $x'\neq\eq$, a contradiction.}
so the notion of \acl{VS} tacitly incorporates a certain degree of maximality.
This is made clearer in the following:

\begin{proposition}
\label{prop:Nash-stable}
Suppose that $\eqset\subseteq\feas$ is variationally stable.
Then:
\begin{enumerate}
[\indent a\upshape)]
\item
$\eqset$ is convex.
\item
If $\game$ is concave, $\eqset$ is an isolated component of \aclp{NE}.
\item
If $\eqset$ is globally stable, it coincides with the game's set of \aclp{NE}.
\end{enumerate}
\end{proposition}


\begin{Proof}[Proof of \cref{prop:Nash-stable}]
To show that $\eqset$ is convex, take $\eq_{0},\eq_{1}\in\eqset$ and set $\eq_{\lambda} = (1-\lambda) \eq_{0} + \lambda \eq_{1}$ for $\lambda\in[0,1]$.
Substituting in \eqref{eq:VS}, we get $\braket{\payv(\eq_{\lambda})}{\eq_{\lambda} - \eq_{0}} = \lambda\braket{\payv(\eq_{\lambda})}{\eq_{1} - \eq_{0}} \leq 0$ and $\braket{\payv(\eq_{\lambda})}{\eq_{\lambda} - \eq_{1}} = -(1-\lambda) \braket{\payv(\eq_{\lambda})}{\eq_{1} - \eq_{0}} \leq 0$, implying that $\braket{\payv(\eq_{\lambda})}{\eq_{1} - \eq_{0}} = 0$.
Writing $\eq_{1} - \eq_{0} = \lambda^{-1} (\eq_{\lambda} - \eq_{0})$, we then get $\braket{\payv(\eq_{\lambda})}{\eq_{\lambda} - \eq_{0}} = 0$.
By \eqref{eq:VS}, we must have $\eq_{\lambda}\in\eqset$ for all $\lambda\in[0,1]$, implying in turn that $\eqset$ is convex.

We now proceed to show that $\eqset$ only consists of \aclp{NE}.
To that end, asssume first that $\eqset$ is globally stable, pick some $\eq\in\eqset$, and let $z_{\play} = x_{\play} - \eq_{\play}$ for some $x_{\play}\in\feas_{\play}$, $\play\in\players$.
Then, for all $\tau\in[0,1]$, we have
\begin{flalign}
\frac{d}{d\tau} \pay_{\play}(\eq_{\play} + \tau z_{\play};\eq_{-\play})
	&= \braket{\payv_{\play}(\eq_{\play} + \tau z_{\play};\eq_{-\play})}{z_{\play}}
	\notag\\
	&= \frac{1}{\tau} \braket{\payv_{\play}(\eq_{\play} +\tau z_{\play};\eq_{-\play})}{\eq_{\play} +\tau z_{\play} - \eq_{\play}}
	\leq 0,
\end{flalign}
where the last inequality follows from \eqref{eq:VS}.
In turn, this shows that $\pay_{\play}(\eq_{\play};\eq_{-\play}) \geq \pay_{\play}(\eq_{\play}+z_{\play};\eq_{-\play}) = \pay_{\play}(x_{\play};\eq_{-\play})$ for all $x_{\play}\in\feas_{\play}$, $\play\in\players$, \ie $\eq$ is a \acl{NE}.
Our claim for locally stable sets then follows by taking $\tau=0$ above and applying \cref{prop:Nash-var}.

We are left to show that there are no other \aclp{NE} close to $\eqset$ (locally or globally).
To do so, assume first that $\eqset$ is locally stable and let $x'\notin\eqset$ be a \acl{NE} lying in a neighborhood $U$ of $\eqset$ where \eqref{eq:VS} holds.
By \cref{prop:Nash-var}, we have $\braket{\payv(x')}{x - x'} \leq 0$ for all $x\in\feas$.
However, since $x'\notin\eqset$, \eqref{eq:VS} implies that $\braket{\payv(x')}{\eq - x'} > 0$ for all $\eq\in\eqset$, a contradiction.
We conclude that there are no other equilibria of $\game$ in $U$, \ie $\eqset$ is an isolated set of \aclp{NE};
the global version of our claim then follows by taking $U=\feas$.
\end{Proof}

\subsection{Tests for \acl{VS}}
\label{sec:test}

We close this section with a second derivative criterion that can be used to verify whether \eqref{eq:VS} holds.
To state it, define the \emph{Hessian} of a game $\game$ as the block matrix $\hessmat(x) = (\hessmat_{\play\playalt}(x))_{\play,\playalt\in\players}$ with
\begin{equation}
\label{eq:Hessian}
\hessmat_{\play\playalt}(x)
	= \tfrac{1}{2} \nabla_{x_{\playalt}} \nabla_{x_{\play}} \pay_{\play}(x)
	+ \tfrac{1}{2} \parens{\nabla_{x_{\play}} \nabla_{x_{\playalt}} \pay_{\playalt}(x)}^{\top}.
\end{equation}
We then have:

\begin{proposition}
\label{prop:Hessian}
If $\eq$ is a \acl{NE} of $\game$ and $\hessmat(\eq)\ml0$ on $\tcone(\eq)$, then $\eq$ is stable \textendash\ and hence an isolated \acl{NE}.
In particular, if $\hessmat(x)\ml0$ on $\tcone(x)$ for all $x\in\feas$, $\eq$ is globally stable \textendash\ so it is the unique equilibrium of $\game$.
\end{proposition}

\begin{remark*}
The requirement ``$\hessmat(\eq)\ml0$ on $\tcone(\eq)$'' above means that $z^{\top} \hessmat(\eq) z < 0$ for every nonzero tangent vector $z\in\tcone(\eq)$.
\end{remark*}

\begin{Proof}
Assume first that $\hessmat(x) \ml 0$ on $\tcone(x)$ for all $x\in\feas$.
By Theorem 6 in \cite{Ros65}, $\game$ satisfies \eqref{eq:MC} so our claim follows from \cref{cor:stable}.
For our second claim, if $\hessmat(\eq) \ml 0$ on $\tcone(\eq)$ for some \acl{NE} $\eq$ of $\game$, we also have $\hessmat(x) \ml 0$ for all $x$ in a neighborhood $U=\prod_{\play\in\players} U_{\play}$ of $\eq$ in $\feas$.
By the same theorem in \cite{Ros65}, we get that \eqref{eq:MC} holds locally in $U$, so the above reasoning shows that $\eq$ is the unique equilibrium of the restricted game $\game\vert_{U}(\players,U,\pay\vert_{U})$.
Hence, $\eq$ is locally stable and isolated in $\game$.
\end{Proof}

We provide two straightforward applications of \cref{prop:Hessian} below:

\begin{example}[Potential games]
Following \cite{MS96}, a game $\game$ is called a \emph{potential game} if it admits a \emph{potential function} $\pot\from\feas\to\R$ such that
\begin{equation*}
\label{eq:potential}
\tag{PF}
\pay_{\play}(x_{\play};x_{-\play}) - \pay_{\play}(x_{\play}';x_{-\play})
	= \pot(x_{\play};x_{-\play}) - \pot(x_{\play}';x_{-\play})
	\quad
	\text{for all $x,x'\in\feas$, $\play\in\players$}.
\end{equation*}
Local maximizers of $\pot$ are \aclp{NE} and the converse also holds if $\pot$ is concave \citep{Ney97}.
By differentiating \eqref{eq:potential}, it is easy to see that the Hessian of $\game$ is just the Hessian of its potential.
Hence, if a game admits a concave potential $\pot$, the game's Nash set $\eqset = \argmax_{x\in\feas} \pot(x)$ is globally stable.
\end{example}

\begin{example}[Cournot revisited]
\label{ex:Cournot-stable}
Consider again the Cournot oligopoly model of \cref{ex:Cournot}.
A simple differentiation yields
\begin{equation}
\label{eq:Cournot-Hessian}
\hessmat_{\play\playalt}(x)
	= \frac{1}{2} \frac{\pd^{2}\pay_{\play}}{\pd x_{\play}\pd x_{\playalt}}
	+ \frac{1}{2} \frac{\pd^{2}\pay_{\playalt}}{\pd x_{\playalt}\pd x_{\play}}
	= -b_{\play} \delta_{\play\playalt}
	- \tfrac{1}{2}(b_{\play} + b_{\playalt}),
\end{equation}
where $\delta_{\play\playalt} = \one\{\play=\playalt\}$ is the Kronecker delta.
This shows that a Cournot oligopoly admits a unique, globally stable equilibrium whenever the RHS of \eqref{eq:Cournot-Hessian} is negative-definite.
This is always the case if the model is symmetric ($b_{\play} = b$ for all $\play\in\players$), but not necessarily otherwise.%
\footnote{This is so because, in the symmetric case, the RHS of \eqref{eq:Cournot-Hessian} is a circulant matrix with eigenvalues $-b$ and $-(\nPlayers+1)b$.}
Quantitatively, if the coefficients $b_{\play}$ are \ac{iid} on $[0,1]$, a Monte Carlo simulation shows that \eqref{eq:Cournot-Hessian} is negative-definite with probability between $65\%$ and $75\%$ for $\nPlayers\in\{2,\dotsc,100\}$.
\end{example}

%% file: Figures/Nash.tex

\colorlet{TangentColor}{DodgerBlue!40!MidnightBlue}
\colorlet{PolarColor}{FireBrick}

\begin{tikzpicture}
[>=stealth,
vecstyle/.style = {->, line width=.5pt},
edgestyle/.style={-, line width=.5pt},
nodestyle/.style={circle, fill=Black,inner sep = .5pt},
plotstyle/.style={color=DarkGreen!80!Cyan,thick}]

\def\radius{2.75}
\def\costhirty{0.8660256}
\def\cosfortyfive{0.7071068}
\def\cosseventyfive{0.2588190451}
\def\diff{0.77886966103678}
\def\veclength{.5}
\def\upangle{125}
\def\downangle{250}
\def\diffangle{45}
\def\midangle{.5*\upangle+.5*\downangle-180}
\def\conescale{.8}

\coordinate (eq) at (0,0);
\coordinate (x) at (-.4*\radius,-.25*\radius);
\coordinate (X) at (-.6*\radius,0.1*\radius);

\fill [TangentColor!5] ++(\downangle:\conescale*\radius) -- (eq.center) --++(\upangle:\conescale*\radius);
\fill [PolarColor!5] ++(\downangle+90:\conescale*\radius) -- (eq.center) --++(\upangle-90:\conescale*\radius);


\draw
[thin, black]
(eq) arc [start angle = \upangle-90, end angle = \upangle-.8*\diffangle, radius = \radius]
node (p1) {};
\draw
[thin, black]
(eq) arc [start angle = \downangle+90, end angle = \downangle+\diffangle, radius = \radius]
node (p2) {};

\fill
[pattern = north west lines, pattern color = black!10]
(p2.center) -- (p1.center) 
arc [start angle = \upangle-.8*\diffangle, end angle = \upangle-90, radius = \radius]
arc [start angle = \downangle+90, end angle = \downangle+\diffangle, radius = \radius]
-- (p2.center);

\draw [thick, black] (eq) arc [start angle = \upangle-90, end angle = \upangle-.8*\diffangle, radius = \radius];
\draw [thick, black] (eq) arc [start angle = \downangle+90, end angle = \downangle+\diffangle, radius = \radius];

\node (X) at (X) {$\feas$};

\draw [edgestyle, thick, TangentColor] (eq.center) -- ++(\upangle:\conescale*\radius);
\draw [edgestyle, thick, TangentColor] (eq.center) -- ++(\downangle:\conescale*\radius) node [left] {$\tcone(\eq)$};

\draw [edgestyle, thick, PolarColor] (eq.center) -- ++(\upangle-90:\conescale*\radius) node [right] {$\pcone(\eq)$};
\draw [edgestyle, thick, PolarColor] (eq.center) -- ++(\downangle+90:\conescale*\radius);

\node [nodestyle, label = above:$\;\;\;\eq$] (eq) at (eq.north) {\phantom{\:}};
\node [nodestyle, label = {[label distance = 3] 90:$\act$}] (x) at (x) {};
\draw [vecstyle] (eq.center) -- (x.center);
\draw [vecstyle] (eq.center) -- (\midangle:1.5) node [above, black] {$\payv(\eq)$};

\end{tikzpicture}

%% file: Tables/Stability.tex

\begin{tabular}{rcll}
\;
	&\;
	&First-order requirement\;
	&Second-order test
	\\
\hline
\acl{NE}
	&\eqref{eq:Nash}
	&$\braket{\payv(\eq)}{x-\eq} \leq 0$
	&N/A
	\\
\hline
Variational stability
	&\eqref{eq:VS}
	&$\braket{\payv(x)}{x-\eq} \leq 0$
	&$\hessmat(\eq) \ml 0$
	\\
\hline
Monotonicity
	&\eqref{eq:MC}
	&$\braket{\payv(x') - \payv(x)}{x'-x} \leq 0$
	\quad
	&$\hessmat(x) \ml 0$
	\\
\hline
Concave potential
	&\eqref{eq:potential}
	&$\payv(x) = \nabla\pot(x)$
	&$\nabla^{2}\pot(x) \ml 0$
	\\
\hline
\end{tabular}

%% file: Learning.tex

\begin{figure}
\centering
\footnotesize
\input{Figures/DA}
\caption{Schematic representation of \acl{DA}.}
\label{fig:DA}
\end{figure}
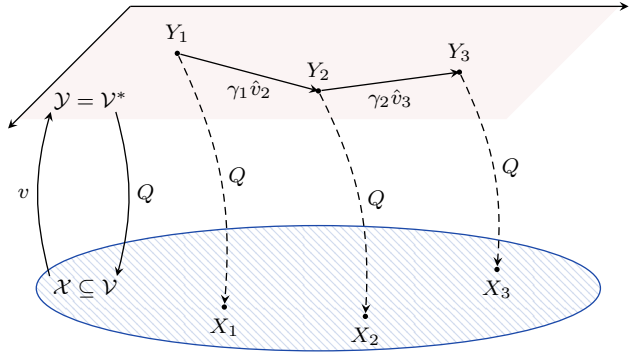

In this section, we adapt the widely used \acdef{DA} method of \cite{Nes09} to our game-theoretic setting.%
\footnote{In optimization, the roots of the method can be traced back to \citef{NY83};
see also \citef{BecTeb03}, \citef{NJLS09} and \citef{SS11}.}
Intuitively, the main idea is as follows:
At each stage of the process, every player $\play\in\players$ gets an estimate $\est_{\play}$ of the individual gradient of their payoff function at the current action profile, possibly subject to noise and uncertainty.
Subsequently, they take a step along this estimate in the dual space $\dspace_{\play}$ (where gradients live),
and they ``mirror'' the output back to the primal space $\feas_{\play}$ in order to choose an action for the next stage and continue playing (for a schematic illustration, see \cref{fig:DA}).

Formally, starting with some arbitrary (and possibly uninformed) gradient estimate $\score_{\start} = \est_{\start}$ at $\run=\start$, this scheme can be described via the recursion
\begin{equation*}
\label{eq:DA}
\tag{DA}
\begin{aligned}
\act_{\play,\run}
	&= \choice_{\play}(\score_{\play,n}),
	\\
\score_{\play,\run+1}
	&= \score_{\play,n} + \step_{\run} \est_{\play,\run+1},
\end{aligned}
\end{equation*}
where:
\begin{enumerate}
[\indent1)]
\item
$\run$ denotes the stage of the process.
\item
$\est_{\play,\run+1} \in \dspace_{\play}$ is an estimate of the individual payoff gradient $\payv_{\play}(\act_{\run})$ of player $\play$ at stage $\run$ (more on this below).
\item
$\score_{\play,\run} \in \dspace_{\play}$ is an auxiliary ``score'' variable that aggregates the $\play$-th player's individual gradient steps.
\item
$\step_{\run}>0$ is a nonincreasing step-size sequence, typically of the form $1/n^{\beta}$ for some $\beta\in(0,1]$.
\item
$\choice_{\play}\from\dspace_{\play}\to\feas_{\play}$ is the \emph{choice map} that outputs the $\play$-th player's action as a function of their score vector $\score_{\play}$ (see below for a rigorous definition).
\end{enumerate}


In view of the above, the core components of \eqref{eq:DA} are
\begin{inparaenum}
[\itshape a\upshape)]
\item
the players' gradient estimates;
and
\item
the choice maps that determine the players' actions.
\end{inparaenum}
In the rest of this section, we discuss both in detail.

\subsection{Feedback and uncertainty}
\label{sec:feedback}

Regarding the players' individual gradient observations, we assume that each player $\play\in\players$ has access to a ``black box'' feedback mechanism \textendash\ an \emph{oracle} \textendash\ which returns an estimate of their payoff gradients at their current action profile.
Of course, this information may be imperfect for a multitude of reasons:
for instance
\begin{inparaenum}
[\itshape i\hspace*{1pt}\upshape)]
\item
estimates may be susceptible to random measurement errors;
\item
the transmission of this information could be subject to noise;
and/or
\item
the game's payoff functions may be stochastic expectations of the form
\begin{equation}
\label{eq:pay-stochastic}
\pay_{\play}(x)	
	= \exof{\hat\pay_{\play}(x;\omega)}
	\quad
	\text{for some random variable $\omega$},
\end{equation}
and players may only be able to observe the realized gradients $\nabla_{x_{\play}} \hat\pay_{\play}(x;\omega)$.
\end{inparaenum}

With all this in mind, we will focus on the noisy feedback model
\begin{equation}
\label{eq:payv-noise}
\est_{\play,\run+1}
	= \payv_{\play}(\act_{\run}) + \noise_{\play,\run+1},
\end{equation}
where the noise process $\noise_{\run} = (\noise_{\play,\run})_{\play\in\players}$ is an $L^{2}$-bounded \acl{MDS} adapted to the history $(\filter_{\run})_{\run=\start}^{\infty}$ of $\act_{\run}$ (\ie $\noise_{\run}$ is $\filter_{\run}$-measurable but $\noise_{\run+1}$ isn't).%
\footnote{Indices have been chosen so that all relevant processes are $\filter_{\run}$-measurable at stage $\run$.}
More explicitly, this means that $\noise_{\run}$ satisfies the statistical hypotheses:
\begin{enumerate}
[\indent 1.]
\item
\emph{Zero-mean:}
\begin{alignat*}{2}
\tag{H1}
\label{eq:zeromean}
\exof*{\noise_{\run+1} \given \filter_{\run}}
	&= 0
	&\quad
	\text{for all $\run=1,2,\dotsc$ \as.}
\intertext{%
\item
\emph{Finite \acl{MSE}:}
there exists some $\noisedev \geq 0$ such that}
\tag{H2}
\label{eq:MSE}
\txs
\exof{\dnorm{\noise_{\run+1}}\crampedllap{^{2}} \given \filter_{\run}}
	&\leq \noisevar
	&\quad
	\text{for all $\run=1,2,\dotsc$ \as.}
\end{alignat*}
\end{enumerate}
Alternatively, \labelcref{eq:zeromean,eq:MSE} simply posit that the players' individual gradient estimates are \emph{conditionally unbiased and bounded in mean square}, viz.
\begin{subequations}
\label{eq:estimates}
\begin{flalign}
\label{eq:unbiased}
&\exof{\est_{\run+1} \given \filter_{\run}}
	= \payv(\act_{\run}),
	\\
\label{eq:vbound}
&\exof{ \dnorm{\est_{\run+1}}^{2} \given \filter_{\run}}
	\leq \vbound^{2}
	\quad
	\text{for some finite $\vbound>0$}.
\end{flalign}
\end{subequations}

The above allows for a broad range of error processes, including all compactly supported, \mbox{(sub-)}Gaussian, \mbox{(sub-)}exponential and log-normal distributions.%
\footnote{In particular, we will not be assuming \acs{iid} errors;
this point is crucial for applications to distributed control where measurements are typically correlated with the state of the system.}
In fact, both hypotheses can be relaxed (for instance, by assuming a small bias or asking for finite moments up to some order $q<2$), but we do not do so to keep things simple.

\subsection{Choosing actions}
\label{sec:choice}

Given that the players' score variables aggregate gradient steps, a reasonable choice for $\choice_{\play}$ would be the $\argmax$ correspondence $y_{\play} \mapsto \argmax_{x_{\play}\in\feas_{\play}} \braket{y_{\play}}{x_{\play}}$ that outputs those actions which are most closely aligned with $y_{\play}$.
Notwithstanding, there are two problems with this approach:
\begin{inparaenum}
[\itshape a\upshape)]
\item
this assignment is too aggressive in the presence of uncertainty;
and
\item
generically, the output would be an extreme point of $\feas$, so \eqref{eq:DA} could never converge to an interior point.
\end{inparaenum}
Thus, instead of taking a ``hard'' $\argmax$ approach, we will focus on \emph{regularized} maps of the form
\begin{equation}
y_{\play}
	\mapsto \argmax_{x_{\play}\in\feas_{\play}} \{ \braket{y_{\play}}{x_{\play}} - h_{\play}(x_{\play}) \},
\end{equation}
where the ``regularization'' term $h_{\play}\from\feas_{\play}\to\R$ satisfies the following requirements:

\begin{definition}
\label{def:choice}
Let $\cvx$ be a compact convex subset of a finite-dimensional normed space $\vecspace$.
We say that $h\from\cvx\to\R$ is a \emph{regularizer} (or \emph{penalty function}) on $\cvx$ if:
\begin{enumerate}
\item
$h$ is continuous.
\item
$h$ is \emph{strongly convex}, \ie there exists some $\strong>0$ such that
\begin{equation}
\label{eq:strong}
h(t x + (1-t) x')
	\leq t h(x) + (1-t) h(x')
	- \tfrac{1}{2} \strong t (1-t) \norm{x' - x}^{2}
\end{equation}
for all $x,x'\in\cvx$ and all $t\in[0,1]$.
\end{enumerate}
The \emph{choice} (or \emph{mirror}) \emph{map $\choice\from\dspace\to\cvx$ induced by $h$} is then defined as
\begin{equation}
\label{eq:choice}
\choice(y)
	= \argmax \setdef{\braket{y}{x} - h(x)}{x\in\cvx}.
\end{equation}
\end{definition}

In what follows, we will be assuming that each player $\play\in\players$ is endowed with an individual penalty function $h_{\play}\from\feas_{\play}\to\R$ that is $\strong_{\play}$-strongly convex.
Furthermore, to emphasize the interplay between primal and dual variables (the players' actions $x_{\play}$ and their score vectors $y_{\play}$ respectively),
we will write $\dual_{\play}\equiv\dspace_{\play}$ for the dual space of $\vecspace_{\play}$ and $\choice_{\play}\from\dual_{\play}\to\feas_{\play}$ for the choice map induced by $h_{\play}$.

More concisely, this information can be encoded in the aggregate penalty function $h(x) = \sum_{\play} h_{\play}(x_{\play})$ with associated strong convexity constant $\strong \equiv \min_{\play} \strong_{\play}$.%
\footnote{We assume here that $\vecspace\equiv\prod_{\play}\vecspace_{\play}$ is endowed with the product norm $\norm{x}_{\vecspace}^{2} = \sum_{\play} \norm{x_{\play}}_{\vecspace_{\play}}^{2}$.}
The induced choice map is simply $\choice \equiv (\choice_{1},\dotsc,\choice_{\nPlayers})$ so we will write $x = \choice(y)$ for the action profile induced by the score vector $y = (y_{1},\dotsc,y_{\nPlayers})\in \dual \equiv \prod_{\play} \dual_{\play}$.

\begin{remark}
In finite games, \citef{MP95} referred to $\choice_{\play}$ as a ``quantal response function'' (the notation $\choice$ alludes precisely to this terminology).
In the same game-theoretic context, the composite map $\choice_{\play}\circ\payv_{\play}$ is often called a smooth, perturbed, or regularized best response;
for a detailed discussion, see \citef{HS02} and \citef{MS16}.
\end{remark}

\smallskip
We discuss below a few examples of this regularization process:

\begin{example}[Euclidean projections]
\label{ex:Eucl}
Let $h(x) = \frac{1}{2} \norm{x}_{2}^{2}$.
Then, $h$ is $1$-strongly convex with respect to $\norm{\cdot}_{2}$ and the corresponding choice map is the closest point projection
\begin{equation}
\label{eq:choice-Eucl}
\Eucl_{\feas}(y)
	\equiv \argmax_{x\in\feas} \braces[\big]{ \braket{y}{x} - \tfrac{1}{2} \norm{x}_{2}^{2} }
	= \argmin_{x\in\feas}\; \norm{y - x}_{2}^{2}.
\end{equation}
The induced learning scheme (cf.~\cref{alg:DA-Eucl}) may thus be viewed as a multi-agent variant of gradient ascent with lazy projections \citefp{Zin03}.
For future reference, note that $h$ is differentiable on $\feas$ and $\Eucl_{\feas}$ is \emph{surjective} (\ie $\im\Eucl_{\feas} = \feas$).
\end{example}

\input{Projection}

\begin{example}[Entropic regularization]
\label{ex:logit}
Motivated by mixed strategy learning in finite games (\cref{ex:finite}), let $\simplex = \setdef{x\in\R^{d}_{+}}{\sum_{j=1}^{d} x_{j} = 1}$ denote the unit simplex of $\R^{d}$.
Then, a standard regularizer on $\simplex$ is provided by the (negative) Gibbs entropy
\begin{equation}
\label{eq:entropy}
h(x)
	= \sum_{\ell=1}^{d} x_{\ell} \log x_{\ell}.
\end{equation}
The entropic regularizer \eqref{eq:entropy} is $1$-strongly convex with respect to the $L^{1}$-norm on $\R^{d}$.
Moreover, a straightforward calculation shows that the induced choice map is
\begin{equation}
\label{eq:choice-logit}
\logit(y)
	= \frac{1}{\sum_{\ell=1}^{d} \exp(y_{\ell})} (\exp(y_{1}),\dotsc,\exp(y_{d})).
\end{equation}
This model is known as \emph{logit choice} and the associated learning scheme has been studied extensively in evolutionary game theory and online learning;
for a detailed account, see \citef{Vov90}, \citef{LW94}, \citef{LM13}, and references therein.
In contrast to the previous example, $h$ is differentiable \emph{only} on the relative interior $\intsimplex$ of $\simplex$ and $\im\logit = \intsimplex$ (\ie $\logit$ is ``essentially'' surjective).
\end{example}

\subsection{Surjectivity vs. steepness}
\label{sec:steepness}

We close this section with an important link between the boundary behavior of penalty functions and the surjectivity of the induced choice maps.
To describe it, it will be convenient to treat $h$ as an extended-real-valued function $h\from\vecspace\to\R\cup\{\infty\}$ by setting $h=\infty$ outside $\feas$.
The \emph{subdifferential} of $h$ at $x\in\vecspace$ is then defined as
\begin{equation}
\label{eq:subdiff}
\pd h(x)
	= \setdef{y\in\dspace}{\text{$h(x') \geq h(x) + \braket{y}{x' - x}$ for all $x'\in\vecspace$}},
\end{equation}
and $h$ is called \emph{subdifferentiable} at $x\in\feas$ whenever $\pd h(x)$ is nonempty.
This is always the case if $x\in\intacts$, so $\intacts \subseteq \dom\pd h \equiv \setdef{x\in\feas}{\pd h(x)\neq\varnothing} \subseteq \feas$ \citefp[Chap.~26]{Roc70}.

Intuitively, $h$ fails to be subdifferentiable at a boundary point $x\in\bd(\feas)$ only if it becomes ``infinitely steep'' near $x$.
We thus say that $h$ is \emph{steep} at $x$ whenever $x\notin\dom h$;
otherwise, $h$ is said to be \emph{nonsteep} at $x$.
The following proposition shows that regularizers that are everywhere nonsteep (as in \cref{ex:Eucl}) induce choice maps that are surjective;
on the other hand, regularizers that are everywhere steep (cf. \cref{ex:logit}) induce choice maps that are interior-valued:

\begin{proposition}
\label{prop:choice}
Let $h$ be a $\strong$-strongly convex regularizer with induced choice map $\choice\from\dual\to\feas$,
and let $h^{\ast}\from\dual\to\R$ be the convex conjugate of $h$, \ie
\begin{equation}
\label{eq:conjugate}
h^{\ast}(y)
	= \max \setdef{\braket{y}{x} - h(x)}{x\in\feas},
	\quad
	\text{$y\in\dual$}.
\end{equation}
Then:
\begin{enumerate}
[\indent\upshape 1)]
\item
$x=\choice(y)$ if and only if $y\in\pd h(x)$;
in particular, $\im\choice = \dom\pd h$.
\item
$h^{\ast}$ is differentiable on $\dual$ and $\nabla h^{\ast}(y) = \choice(y)$ for all $y\in\dual$.
\item
$\choice$ is $(1/\strong)$-Lipschitz continuous.
\end{enumerate}
\end{proposition}

\cref{prop:choice} is essentially folklore in optimization and convex analysis;
for a proof, see \citet[Theorem 23.5]{Roc70} and \citet[Theorem 12.60(b)]{RW98}.

%% file: Figures/DA.tex

\colorlet{TangentColor}{DodgerBlue!40!MidnightBlue}
\colorlet{PolarColor}{FireBrick}

\begin{tikzpicture}
[auto,
>=stealth,
vecstyle/.style = {->, line width=.5pt},
edgestyle/.style={-, line width=.5pt},
nodestyle/.style = {circle,fill=black,inner sep=0, minimum size=2},
plotstyle/.style={color=DarkGreen!80!Cyan,thick}]

\def\radius{2.5}
\def\costhirty{0.8660256}
\def\cosfortyfive{0.7071068}
\def\diff{0.77886966103678}
\def\veclength{.5}
\def\conescale{.8}

\coordinate (base) at (0,0);
\coordinate (dbase) at (0,1.5*\radius);

\coordinate (X) at ($(base) + (-.25*\radius,-0*\radius)$);
\coordinate (Y) at ($(dbase) + (-.25*\radius,-.5*\radius)$);

\coordinate (y1) at ($(dbase) + (.25*\radius,-.25*\radius)$);
\coordinate (y2) at ($(dbase) + (1*\radius,-.45*\radius)$);
\coordinate (y3) at ($(dbase) + (1.75*\radius,-.35*\radius)$);
\coordinate (y4) at ($(dbase) + (2.25*\radius,-.1*\radius)$);

\coordinate (x1) at ($(base) + (.5*\radius,-.1*\radius)$);
\coordinate (x2) at ($(base) + (1.25*\radius,-.15*\radius)$);
\coordinate (x3) at ($(base) + (1.95*\radius,.1*\radius)$);
\coordinate (x4) at ($(base) + (2.25*\radius,0.1*\radius)$);

\filldraw
[pattern = north west lines, pattern color = TangentColor!20, draw=TangentColor, edgestyle]
($(base) + (\radius,0)$) ellipse [x radius= 1.5*\radius, y radius= \radius/3];

\fill
[PolarColor!5]
(dbase) -- ($(dbase) - (.6*\radius,.6*\radius)$) -- ($(dbase) + (2*\radius,-.6*\radius)$)-- ($(dbase) + (2.6*\radius,0)$);
\draw
[vecstyle]
(dbase) -- ($(dbase) + (2.65*\radius,0)$);
\draw
[vecstyle]
(dbase) -- ($(dbase) - .65*(\radius,\radius)$);

\node [inner sep = 1pt] (X) at (X) {$\feas \subseteq \vecspace\!$};
\node [inner sep = 1pt] (Y) at (Y) {$\dual = \dspace\!\!\!\!$};
\draw [vecstyle] (X.north west) to [bend left =15] node [left] {$\payv$}(Y.south west);
\draw [vecstyle] (Y.south east) to [bend left =15] node [right] {$\choice$}(X.north east);

\node [nodestyle, label = above:$\score_{1}$] (y1) at (y1) {};
\node [nodestyle, label = above:$\score_{2}$] (y2) at (y2) {};
\node [nodestyle, label = above:$\score_{3}$] (y3) at (y3) {};

\draw [vecstyle] (y1) -- (y2.center) node [below,midway,black] {$\step_{1} \est_{2}$};
\draw [vecstyle] (y2) -- (y3.center) node [below,midway,black] {$\step_{2} \est_{3}$};

\node [nodestyle, label = below:$\act_{1}$] (x1) at (x1) {};
\node [nodestyle, label = below:$\act_{2}$] (x2) at (x2) {};
\node [nodestyle, label = below:$\act_{3}$] (x3) at (x3) {};

\draw [vecstyle, densely dashed] (y1) to [bend left = 15] node [right] {$\choice$} (x1);
\draw [vecstyle, densely dashed] (y2) to [bend left = 15] node [right] {$\choice$} (x2);
\draw [vecstyle, densely dashed] (y3) to [bend left = 15] node [right] {$\choice$} (x3);


\end{tikzpicture}

%% file: Projection.tex

\begin{algorithm}[tbp]
\caption{Dual averaging with Euclidean projections (\cref{ex:Eucl}).}
\label{alg:DA-Eucl}
\normalsize
\begin{algorithmic}[1]
	\REQUIRE
	step-size sequence $\step_{\run}\propto 1/\run^{\beta}$, $\beta\in(0,1]$;
	initial scores $\score_{\play}\in\dual_{\play}$
	\FOR{$\run=\start,2,\dotsc$}
		\FOR{every player $\play\in\players$}
			\STATE play $\act_{\play} \leftarrow \Eucl_{\feas_{\play}}(\score_{\play})$;
			\hfill
			\COMMENT{choose an action}
			\\[2pt]
			\STATE observe $\est_{\play}$;
			\hfill
			\COMMENT{estimate gradient}
			\\[2pt]
			\STATE update $\score_{\play} \leftarrow \score_{\play} + \step_{\run} \est_{\play}$;
			\hfill
			\COMMENT{take gradient step}
		\ENDFOR
	\ENDFOR
\end{algorithmic}
\end{algorithm}

%% file: Analysis.tex

A key property of \eqref{eq:DA} in concave games is that it leads to \emph{no regret}, viz.
\begin{equation}
\label{eq:regret}
\max_{x_{\play}\in\feas_{\play}}
	\sum_{\iRun=\start}^{\run} \left[ \pay_{\play}(x_{\play};\act_{-\play,\iRun}) - \pay_{\play}(\act_{\iRun}) \right]
	= o(\run)
	\quad
	\text{for all $\play\in\players$},
\end{equation}
provided that
the algorithm's step-size is chosen appropriately \textendash\ for a precise statement, see \cite{Xia10} and \cite{SS11}.
As such, under \eqref{eq:DA}, every player's average payoff matches asymptotically that of the best fixed action in hindsight (though, of course, this does not take into account changes to other players' actions due to a change in a given player's chosen action).

In this section, we expand on this worst-case guarantee and we derive some general convergence results for the actual sequence of play induced by \eqref{eq:DA}.
Specifically, in \cref{sec:limits} we show that if \eqref{eq:DA} converges to some action profile, this limit is a \acl{NE}.
Subsequently, to obtain stronger convergence results, we introduce in \cref{sec:Fenchel} the so-called \emph{Fenchel coupling}, a ``primal-dual'' divergence measure between the players' (primal) action variables $x_{\play}\in\feas_{\play}$ and their (dual) score vectors $y_{\play}\in\dual_{\play}$.
Using this coupling as a Lyapunov function, we show in \cref{sec:global,sec:local} that globally (resp. locally) stable states are globally (resp. locally) attracting under \eqref{eq:DA}.
Finally, in \cref{sec:zerosum}, we examine the convergence properties of \eqref{eq:DA} in zero-sum concave-convex games.

\subsection{Limit states}
\label{sec:limits}

We first show that if the sequence of play induced by \eqref{eq:DA} converges to some $\eq\in\feas$ with positive probability, this limit is a \acl{NE}:

\begin{theorem}
\label{thm:terminal}
Suppose that \eqref{eq:DA} is run with
imperfect gradient information satisfying \crefrange{eq:zeromean}{eq:MSE}
and
a step-size sequence $\step_{\run}$ such that
\begin{equation}
\label{eq:step-loose}
\sum_{\run=\start}^{\infty}
	\Big( \frac{\step_{\run}}{\tau_{\run}} \Big)^{2}
	< \sum_{\run=\start}^{\infty} \step_{\run}
	= \infty,
\end{equation}
where $\tau_{\run} = \sum_{\iRun=\start}^{\run} \step_{\iRun}$.
If the game is \textpar{pseudo-}concave and $\act_{\run}$ converges to $\eq\in\feas$ with positive probability, $\eq$ is a \acl{NE}.
\end{theorem}

\begin{remark}
Note here that the requirement \eqref{eq:step-loose} holds for every step-size policy of the form $\step_{\run} \propto 1/n^{\beta}$, $\beta\leq1$ (i.e. even for \emph{increasing} $\step_{\run}$).
\end{remark}

\begin{Proof}[Proof of \cref{thm:terminal}]
Let $\payveq = \payv(\eq)$ and assume ad absurdum that $\eq$ is not a \acl{NE}.
By the characterization \eqref{eq:Nash-vec} of \aclp{NE}, there exists a player $\play\in\players$ and a deviation $\dev_{\play}\in\feas_{\play}$ such that $\braket{\payveq_{\play}}{\dev_{\play} - \eq_{\play}} > 0$.
Thus, by continuity, there exists some $a>0$ and neighborhoods $U$, $V$ of $\eq$ and $\payveq$ respectively, such that
\begin{equation}
\label{eq:notNash}
\braket{\payv_{\play}'}{\dev_{\play} - x_{\play}'}
	\geq \sharp
\end{equation}
whenever $x'\in U$ and $\payv'\in V$.

Now, let $\Omega_{0}$ be the event that $\act_{\run}$ converges to $\eq$, so $\probof{\Omega_{0}} > 0$ by assumption.
Within $\Omega_{0}$, we may assume for simplicity that $\act_{\run}\in U$ and $\payv(\act_{\run})\in V$ for all $\run$, so \eqref{eq:DA} yields
\begin{flalign}
\score_{\run+1}
	&= \score_{\start} + \sum_{\iRun=\start}^{\run} \step_{\iRun} \est_{\iRun+1}
	= \score_{\run_{0}} + \sum_{\iRun=\start}^{\run} \step_{\iRun} \left[ \payv(\act_{\iRun}) + \noise_{\iRun+1} \right]
	= \score_{\start} + \tau_{\run} \bar\payv_{\run+1},
\end{flalign}
where we set $\bar\payv_{\run+1} = \tau_{\run}^{-1} \sum_{\iRun=\start}^{\run} \step_{\iRun} \est_{\iRun+1} = \tau_{\run}^{-1} \sum_{\iRun=\start}^{\run} \step_{\iRun} \left[ \payv(\act_{\iRun}) +\noise_{\iRun+1} \right]$.

We now claim that $\probof{\bar\payv_{\run} \to \payveq \given \Omega_{0}} = 1$.
Indeed, by \eqref{eq:step-loose} and \eqref{eq:MSE}, we have
\begin{equation}
\sum_{\run=\start}^{\infty} \frac{1}{\tau_{\run}^{2}} \exof{ \dnorm{\step_{\run} \noise_{\run+1}}^{2} \given \filter_{\run}}
	\leq \sum_{\run=\start}^{\infty} \frac{\step_{\run}^{2}}{\tau_{\run}^{2}} \noisevar
	< \infty.
\end{equation}
Therefore, by the \acl{LLN} for \aclp{MDS} \citep[Theorem 2.18]{HH80}, we obtain $\tau_{\run}^{-1}\sum_{\iRun=\start}^{\run} \step_{\iRun} \noise_{\iRun+1} \to 0$ \as.
Given that $\payv(\act_{\run})\to\payveq$ in $\Omega_{0}$ and $\probof{\Omega_{0}} > 0$, we infer that $\probof{\bar\payv_{\run} \to \payveq \given \Omega_{0}} = 1$, as claimed.

Now, with $\score_{\play,\run} \in \pd h_{\play}(\act_{\play,\run})$ by \cref{prop:choice}, we also have
\begin{equation}
\label{eq:hdiff}
h_{\play}(\dev_{\play}) - h_{\play}(\act_{\play,\run})
	\geq \braket{\score_{\play,\run}}{\dev_{\play} - \act_{\play,\run}}
	= \braket{\score_{\play,\start}}{\dev_{\play} - \act_{\play,\run}} + \tau_{\run-1} \braket{\bar\payv_{\play,\run}}{\dev_{\play} - \act_{\play,\run}}.
\end{equation}
Since $\bar\payv_{\run}\to\payveq$ almost surely on $\Omega_{0}$, \eqref{eq:notNash} yields $\braket{\bar\payv_{\play,\run}}{\dev_{\play} - \act_{\play,\run}} \geq \sharp > 0$ for all sufficiently large $\run$.
However,
given that $\abs{\braket{\score_{\play,\start}}{\dev_{\play} - \act_{\play,\run}}} \leq \dnorm{\score_{\play,\start}} \norm{\dev_{\play} - \act_{\play,\run}} \leq \dnorm{\score_{\play,\start}} \norm{\feas} = \bigoh(1)$,
a simple substitution in \eqref{eq:hdiff} yields $h_{\play}(\dev_{\play}) - h_{\play}(\act_{i,\run}) \gtrsim \sharp \tau_{\run}\to\infty$ with positive probability, a contradiction.
We conclude that $\eq$ is a \acl{NE} of $\game$, as claimed.
\end{Proof}

\subsection{The Fenchel coupling}
\label{sec:Fenchel}

A key tool in establishing the convergence properties of \eqref{eq:DA} is the so-called \emph{Bregman divergence} $\breg(\base,\notbase)$ between a given base point $\base\in\feas$ and a test state $\notbase\in\feas$.
Following \cite{Kiw97b}, $\breg(\base,\notbase)$ is defined as the difference between $h(\base)$ and the best linear approximation of $h(\base)$ from $x$, viz.
\begin{equation}
\label{eq:Bregman}
\breg(\base,\notbase)
	= h(\base) - h(\notbase) - h'(\notbase;\base-\notbase),
\end{equation}
where $h'(\notbase;z) = \lim_{t\to0^{+}} t^{-1} [h(\notbase+tz) - h(\notbase)]$ denotes the one-sided derivative of $h$ at $\notbase$ along $z\in\tcone(\notbase)$.
Owing to the (strict) convexity of $h$, we have $\breg(\base,\notbase) \geq 0$ and $\act_{\run}\to\base$ whenever $\breg(\base,\act_{\run})\to0$ \citep{Kiw97b}.
Accordingly, the convergence of a sequence $\act_{\run}$ to a target point $\base$ can be checked directly by means of the associated divergence $\breg(\base,\act_{\run})$.

Nevertheless, it is often impossible to glean any useful information on $\breg(\base,\act_{\run})$ from \eqref{eq:DA} when $\act_{\run} = \choice(\score_{\run})$ is not interior.
Instead, given that \eqref{eq:DA} mixes primal and dual variables (actions and scores respectively),
it will be more convenient to use the following ``primal-dual divergence'' between dual vectors $y\in\dual$ and base points $\base\in\feas$:

\begin{definition}
\label{def:Fenchel}
Let $h\from\feas\to\R$ be a penalty function on $\feas$.
Then, the \emph{Fenchel coupling} induced by $h$ is defined as
\begin{equation}
\label{eq:Fenchel}
\fench(\base,y)
	= h(\base) + h^{\ast}(y) - \braket{y}{\base}
	\quad
	\text{for all $\base\in\feas$, $y\in\dual$}.
\end{equation}
\end{definition}

The terminology ``Fenchel coupling'' is due to \cite{MS16} and refers to the fact that \eqref{eq:Fenchel} collects all terms of Fenchel's inequality.
As a result, $\fench(\base,y)$ is nonnegative and strictly convex in both arguments (though not jointly so).
Moreover, it enjoys the following key properties:

\begin{proposition}
\label{prop:Fenchel}
Let $h$ be a $\strong$-strongly convex penalty function on $\feas$.
Then, for all $\base\in\feas$ and all $y,y'\in\dual$, we have:
\begin{subequations}
\label{eq:Fench-properties}
\begin{alignat}{2}
\label{eq:Fench-Bregman}
&a)
	\;\;
	\fench(\base,y)
	&&= \breg(\base,\choice(y))
	\;\;
	\text{if $\choice(y)\in\intacts$ \textpar{but not necessarily otherwise}}.
	\\[2pt]
\label{eq:Fench-norm}
&b)
	\;\;
	\fench(\base,y)
	&&\geq \tfrac{1}{2} \strong \, \norm{\choice(y) - \base}^{2}.
	\\[2pt]
\label{eq:Fench-bound}
&c)
	\;\;
	\fench(\base,y')
	&&\leq \fench(\base,y) + \braket{y' - y}{\choice(y) - \base} + \tfrac{1}{2\strong} \dnorm{y'-y}^{2}.
\end{alignat}
\end{subequations}
\end{proposition}

\cref{prop:Fenchel} (proven in \cref{app:aux}) justifies the terminology ``primal-dual divergence'' and plays a key role in our analysis.
Specifically, given a sequence $\score_{\run}$ in $\dual$, \eqref{eq:Fench-norm} yields $\choice(\score_{\run}) \to \base$ whenever $\fench(\base,\score_{\run}) \to 0$, meaning that $\fench(\base,\score_{\run})$ can be used to test the convergence of $\choice(\score_{\run})$ to $\base$.

For technical reasons, it is convenient to also assume the converse, namely that
\begin{equation}
\label{eq:Fench-reg}
\tag{H3}
\fench(\base,\score_{\run}) \to 0
	\quad
	\text{whenever}
	\quad
	\choice(\score_{\run}) \to \base.
\end{equation}
When $h$ is steep, we have $\fench(\base,y) = \breg(\base,\choice(y))$ for all $y\in\dual$, so \eqref{eq:Fench-reg} boils down to the requirement
\begin{equation}
\label{eq:Breg-conv}
\breg(\base,\notbase_{\run}) \to 0
	\quad
	\text{whenever}
	\quad
	\act_{\run} \to \base.
\end{equation}
This so-called ``reciprocity condition'' is well known in the theory of Bregman functions \citep{CT93,Kiw97b,ABB04}:
essentially, it means that the sublevel sets of $\breg(\base,\cdot)$ are neighborhoods of $\base$ in $\feas$.
Hypothesis \eqref{eq:Fench-reg} posits that the \emph{images} of the sublevel sets of $\fench(\base,\cdot)$ under $\choice$ are neighborhoods of $\base$ in $\feas$, so it may be seen as a ``primal-dual'' variant of Bregman reciprocity.
Under this light, it is easy to check that \cref{ex:Eucl,ex:logit} both satisfy \eqref{eq:Fench-reg}.

Obviously, when \eqref{eq:Fench-reg} holds, \cref{prop:Fenchel} gives:

\begin{corollary}
Under \eqref{eq:Fench-reg}, $\fench(\base,\score_{\run})\to 0 $ if and only if $\choice(\score_{\run})\to\base$.
\end{corollary}

To extend the above to subsets of $\feas$, we further define the setwise coupling
\begin{equation}
\label{eq:Fench-set}
\fench(\baseset,y)
	= \inf\setdef{\fench(\base,y)}{\base\in\baseset},
	\quad
	\text{$\baseset\subseteq\feas$, $y\in\dual$}.
\end{equation}
In analogy to the pointwise case, we then have:

\begin{proposition}
\label{prop:Fench-set}
Let $\baseset$ be a closed subset of $\feas$.
Then, $\choice(\score_{\run})\to\baseset$ whenever $\fench(\baseset,\score_{\run})\to0$;
in addition, if \eqref{eq:Fench-reg} holds, the converse is also true.
\end{proposition}

The proof of \cref{prop:Fench-set} is a straightforward exercise in point-set topology so we omit it.
What's more important is that, thanks to \cref{prop:Fench-set}, the Fenchel coupling can also be used to test for convergence to a set;
in what follows, we employ this property freely.

%

\subsection{Global convergence}
\label{sec:global}

In this section, we focus on globally stable \aclp{NE} (and sets thereof).
We begin with the perfect feedback case:

\begin{theorem}
\label{thm:global-perfect}
Suppose that \eqref{eq:DA} is run with
perfect feedback \textpar{$\noisedev=0$},
choice maps satisfying \eqref{eq:Fench-reg},
and
a step-size $\step_{\run}$ such that
$\sum_{\iRun=\start}^{\run} \step_{\iRun}^{2} \big/ \sum_{\iRun=\start}^{\run} \step_{\iRun} \to 0$.
If the set $\eqset$ of the game's \aclp{NE} is globally stable, $\act_{\run}$ converges to $\eqset$.
\end{theorem}

\begin{Proof}
Let $\eqset$ be the game's set of \aclp{NE},
fix some arbitrary $\eps>0$,
and let $U_{\eps} = \setdef{x = \choice(y)}{\fench(\eqset,y) < \eps}$.
Then, by \cref{prop:Fench-set}, it suffices to show that $\act_{\run}\in U_{\eps}$ for all sufficiently large $\run$.

To that end,
for all $\eq\in\eqset$, \cref{prop:Fenchel} yields
\begin{equation}
\label{eq:Fbound-det}
\fench(\eq,\score_{\run+1})
	\leq \fench(\eq,\score_{\run}) + \step_{\run} \braket{\payv(\act_{\run})}{\act_{\run} - \eq} + \frac{\step_{\run}^{2}}{2\strong} \dnorm{\payv(\act_{\run})}^{2}.
\end{equation}
To proceed, assume inductively that $\act_{\run}\in U_{\eps}$.
By \eqref{eq:Fench-reg}, there exists some $\delta>0$ such that $\cl(U_{\eps/2})$ contains a $\delta$-neighborhood of $\eqset$.%
\footnote{Indeed, if this were not the case, there would exist a sequence $\score_{\run}'$ in $\dual$ such that $\choice(\score_{\run}')\to\eqset$ but $\fench(\eqset,\score_{\run}')\geq \eps/2$, in contradiction to \eqref{eq:Fench-reg}.}
Consequently, with $\eqset$ globally stable, there exists some $\sharp\equiv \sharp(\eps) > 0$ such that
\begin{equation}
\label{eq:gradbound-U}
\braket{\payv(x)}{x - \eq}
	\leq - \sharp
	\quad
	\text{for all $x\in U_{\eps}-U_{\eps/2}$, $\eq\in\eqset$}.
\end{equation}
If $\act_{\run} \in U_{\eps} - U_{\eps/2}$ and $\step_{\run} \leq 2\sharp\strong/\vbound^{2}$, \eqref{eq:Fbound-det} yields $\fench(\eq,\score_{\run+1}) \leq \fench(\eq,\score_{\run})$.%
\footnote{Since $\noisedev=0$, we can take here $\vbound = \max_{x\in\feas} \dnorm{\payv(x)}$.}
Hence, minimizing over $\eq\in\eqset$, we get $\fench(\eqset,\score_{\run+1}) \leq \fench(\eqset,\score_{\run}) < \eps$, so $\act_{\run+1} = \choice(\score_{\run+1}) \in U_{\eps}$.
Otherwise, if $\act_{\run} \in U_{\eps/2}$ and $\step_{\run}^{2} < \eps\strong/\vbound^{2}$, combining \eqref{eq:VS} with \eqref{eq:Fbound-det} yields $\fench(\eq,\score_{\run+1}) \leq \fench(\eq,\score_{\run}) + \eps/2$ so, again, $\fench(\eqset,\score_{\run+1}) \leq \fench(\eqset,\score_{\run}) + \eps/2 \leq \eps$, i.e. $\act_{\run+1} \in U_{\eps}$.
We thus conclude that $\act_{\run+1}\in U_{\eps}$ whenever $\act_{\run}\in U_{\eps}$ and $\step_{\run} < \min\{2\sharp\strong/\vbound^{2}, \sqrt{\eps\strong}/\vbound\}$.

To complete the proof, \cref{lem:visit} shows that $\act_{\run}$ visits $U_{\eps}$ infinitely often under the stated assumptions.
Since $\step_{\run}\to0$, our assertion follows.
\end{Proof}


\begin{table}[tbp]
\centering
\renewcommand{\arraystretch}{1.3}
\small
\input{Tables/Hypotheses}
\vspace{2ex}
\caption{Overview of the various regularity hypotheses used in the paper.}
\label{tab:rates}
\vspace{-3ex}
\end{table}


We next show that \cref{thm:global-perfect} extends to the case of imperfect feedback under the additional regularity requirement:
\begin{equation}
\label{eq:Lipschitz}
\tag{H4}
\text{The gradient field $\payv(x)$ is Lipschitz continuous.}
\end{equation}
With this extra assumption, we have:

\begin{theorem}
\label{thm:global-imperfect}
Suppose that \eqref{eq:DA} is run with a step-size sequence $\step_{\run}$ such that $\sum_{\run=\start}^{\infty} \step_{\run}^{2} < \infty$ and $\sum_{\run=\start}^{\infty} \step_{\run} = \infty$.
If \crefrange{eq:zeromean}{eq:Lipschitz} hold and
the set $\eqset$ of the game's \aclp{NE} is globally stable, $\act_{\run}$ converges to $\eqset$ \as.
\end{theorem}

\begin{corollary}
\label{cor:monotone}
If $\game$ satisfies \eqref{eq:MC},
$\act_{\run}$ converges to the \textpar{necessarily unique} \acl{NE} of $\game$ \as.
\end{corollary}

\begin{corollary}
\label{cor:potential}
If $\game$ admits a concave potential,
$\act_{\run}$ converges to the set of \aclp{NE} of $\game$ \as.
\end{corollary}

Because of the noise affecting the players' gradient estimates, our proof strategy for \cref{thm:global-imperfect} is quite different from that of \cref{thm:global-perfect}.
In particular, instead of working directly in discrete time, we start with the continuous-time system
\begin{equation}
\label{eq:DA-cont}
\tag{\ref*{eq:DA}-c}
\begin{aligned}
\dot y
	&= \payv(x),
	\\
x
	&= \choice(y),
\end{aligned}
\end{equation}
which can be seen as a ``mean-field'' approximation of the recursive scheme \eqref{eq:DA}.
As we show in \cref{app:aux}, the orbits $x(t) = \choice(y(t))$ of \eqref{eq:DA-cont} converge to $\eqset$ in a certain, ``uniform'' way.
Moreover, under the assumptions of \cref{thm:global-imperfect}, the sequence $\score_{\run}$ generated by the discrete-time, stochastic process \eqref{eq:DA} comprises an \acdef{APT} of the dynamics \eqref{eq:DA-cont}, i.e. $\score_{\run}$ asymptotically tracks the flow of \eqref{eq:DA-cont} with arbitrary accuracy over windows of arbitrary length \cite{Ben99}.%
\footnote{For a precise definition, see \eqref{eq:APT} below.}
\acp{APT} have the key property that, in the presence of a global attractor, they cannot stray too far from the flow of \eqref{eq:DA-cont};
however, given that $\choice$ may fail to be invertible, the trajectories $x(t) = \choice(y(t))$ do \emph{not} consitute a semiflow, so it is not possible to leverage the general \acl{SA} theory of \cite{Ben99}.
To overcome this difficulty, we exploit the derived convergence bound for $x(t) = \choice(y(t))$, and we then use an inductive shadowing argument to show that \eqref{eq:DA} converges itself to $\eqset$.

\begin{Proof}[Proof of \cref{thm:global-imperfect}]
Fix some $\eps>0$,
let $U_{\eps} = \setdef{x=\choice(y)}{\fench(\eqset,y) < \eps}$,
and
write $\semiflow_{t}\from\dual\to\dual$ for the semiflow induced by \eqref{eq:DA-cont} on $\dual$ \textendash\ i.e. $(\semiflow_{t}(y))_{t\geq0}$ is the solution orbit of \eqref{eq:DA-cont} that starts at $y\in\dual$.%
\footnote{That such a trajectory exists and is unique is a consequence of \eqref{eq:Lipschitz}.}

We first claim there exists some finite $\tau\equiv \tau(\eps)$ such that $\fench(\eqset,\semiflow_{\tau}(y)) \leq \max\{\eps,\fench(\eqset,y) - \eps\}$ for all $y\in\dual$.
Indeed, since $\cl(U_{\eps})$ is a closed neighborhood of $\eqset$ by \eqref{eq:Fench-reg},
\eqref{eq:VS} implies that there exists some $\sharp \equiv \sharp(\eps) > 0$ such that
\begin{equation}
\braket{\payv(x)}{x - \eq}
	\leq -\sharp
	\quad
	\text{for all $\eq\in\eqset$, $x\notin U_{\eps}$}.
\end{equation}
Consequently, if $\tau_{y} = \inf\setdef{t>0}{\choice(\semiflow_{t}(y)) \in U_{\eps}}$ denotes the first time at which an orbit of \eqref{eq:DA-cont} reaches $U_{\eps}$, \cref{lem:Lyapunov} in \cref{app:aux} gives:
\begin{equation}
\label{eq:Fench-bound1}
\fench(\eq,\semiflow_{t}(y))
	\leq \fench(\eq,y) - \sharp t
	\quad
	\text{for all $\eq\in\eqset$, $t\leq\tau_{y}$}.
\end{equation}
In view of this, set $\tau=\eps/\sharp$ and consider the following two cases:
\begin{enumerate}
\item
$\tau_{y} \geq \tau$:
then, \eqref{eq:Fench-bound1} gives $\fench(\eq,\semiflow_{\tau}(y)) \leq \fench(\eq,y) - \eps$ for all $\eq\in\eqset$, so $\fench(\eqset,\semiflow_{\tau}(y)) \leq \fench(\eqset,y) - \eps$.
\item
$\tau_{y} < \tau$:
then, $\choice(\semiflow_{\tau}(y)) \in U_{\eps}$, so $\fench(\eqset,\semiflow_{\tau}(y)) \leq \eps$.
\end{enumerate}
In both cases we have $\fench(\eqset,\semiflow_{\tau}(y)) \leq \max\{\eps,\fench(\eqset,y) - \eps\}$, as claimed.

Now, let $(Y(t))_{t\geq0}$ denote the affine interpolation of the sequence $\score_{\run}$ generated by \eqref{eq:DA},
i.e. $Y$ is the continuous curve which joins the values $\score_{\run}$ at all times $\tau_{\run} = \sum_{\iRun=\start}^{\run} \step_{\iRun}$.
Under the stated assumptions, a standard result of \citet[Propositions 4.1 and 4.2]{Ben99} shows that $Y(t)$ is an \acl{APT} of $\semiflow$, i.e.
\begin{equation}
\label{eq:APT}
\lim_{t\to\infty} \sup_{0\leq h \leq T} \dnorm{Y(t+h) - \semiflow_{h}(Y(t))}
	= 0
	\quad
	\text{for all $T>0$ \as}.
\end{equation}
Thus, with some hindsight, let $\delta \equiv \delta(\eps)$ be such that $\delta\norm{\feas} + \delta^{2}/(2\strong) \leq \eps$
and
choose $t_{0} \equiv t_{0}(\eps)$ so that $\sup_{0\leq h\leq \tau} \dnorm{Y(t+h) - \semiflow_{h}(Y(t))} \leq \delta$ for all $t\geq t_{0}$.
Then, for all $t\geq t_{0}$ and all $\eq\in\eqset$, \cref{prop:Fenchel} gives
\begin{flalign}
\fench(\eq,Y(t+h))
	&\leq \fench(\eq,\semiflow_{h}(Y(t)))
	\notag\\
	&+ \braket{Y(t+h) - \semiflow_{h}(Y(t))}{\choice(\semiflow_{h}(Y(t))) - \eq}
	\notag\\
	&+ \frac{1}{2\strong} \dnorm{Y(t+h) - \semiflow_{h}(Y(t))}^{2}
	\notag\\
	&\leq \fench(\eq,\semiflow_{h}(Y(t)))
	+ \delta \norm{\feas}
	+ \frac{\delta^{2}}{2\strong}
	\notag\\
	&\leq \fench(\eq,\semiflow_{h}(Y(t))) + \eps.
\end{flalign}
Hence, minimizing over $\eq\in\eqset$, we get
\begin{equation}
\label{eq:Fench-3eps}
\fench(\eqset,Y(t+h))
	\leq \fench(\eqset,\semiflow_{h}(Y(t))) + \eps
	\quad
	\text{for all $t\geq t_{0}$}.
\end{equation}

By \cref{lem:visit}, there exists some $t\geq t_{0}$ such that $\fench(\eqset,Y(t)) \leq 2\eps$ \as.
Thus, given that $\fench(\eqset,\semiflow_{h}(Y(t)))$ is nonincreasing in $h$ by \cref{lem:Lyapunov}, \cref{eq:Fench-3eps} yields $\fench(\eqset,Y(t+h)) \leq 2\eps + \eps = 3\eps$ for all $h\in[0,\tau]$.
However, by the definition of $\tau$, we also have $\fench(\eqset,\semiflow_{\tau}(Y(t))) \leq \max\{\eps,\fench(\eqset,Y(t)) - \eps\}\leq \eps$, implying in turn that $\fench(\eqset,Y(t+\tau)) \leq \fench(\eqset,\semiflow_{\tau}(Y(t))) + \eps \leq 2\eps$.
Therefore, by repeating the above argument at $t+\tau$ and proceeding inductively, we get $\fench(\eqset,Y(t+h)) \leq 3\eps$ for all $h\in[\iRun\tau,(\iRun+1)\tau]$, $\iRun=1,2,\dotsc$ \as.
Since $\eps$ has been chosen arbitrarily, we conclude that $\fench(\eqset,\score_{\run}) \to 0$, so $\act_{\run}\to\eqset$ by \cref{prop:Fench-set}.
\end{Proof}

We close this section with a few remarks:

\begin{remark}
In the above, the Lipschitz continuity assumption \eqref{eq:Lipschitz} is used to show that the sequence $\act_{\run}$ comprises an \ac{APT} of the continuous-time dynamics \eqref{eq:DA-cont}.
Since any continuous functions on a compact set is uniformly continuous, the proof of Proposition 4.1 in \citet[p.~14]{Ben99} shows that \eqref{eq:Lipschitz} can be dropped altogether if \eqref{eq:DA-cont} is well-posed (which, in turn, holds if $\payv(x)$ is only \emph{locally} Lipschitz).
Albeit less general, Lipschitz continuity is more straightforward as an assumption, so we do not go into the details of this relaxation.

We should also note that several classic convergence results for \acl{DA} and \acl{MD} do not require Lipschitz continuity at all (see e.g. \citealp{Nes09}, and \citealp{NJLS09}).
The reason for this is that these results focus on the convergence of the averaged sequence $\bar \act_{\run} = \sum_{\iRun=\start}^{\run} \step_{\iRun} \act_{\iRun} \big/ \sum_{\iRun=\start}^{\run} \step_{\iRun}$, whereas the figure of merit here is the \emph{actual} sequence of play $\act_{\run}$.
The latter sequence is more sensitive to noise, hence the need for additional regularity;
in our ergodic analysis later in the paper, \eqref{eq:Lipschitz} is not invoked.
\end{remark}

\begin{remark}
\cref{thm:global-imperfect} shows that \eqref{eq:DA} converges to equilibrium, but the summability requirement $\sum_{\run=\start}^{\infty} \step_{\run}^{2} < \infty$ suggests that players must be more conservative under uncertainty.
To make this more precise, note that the step-size assumptions of \cref{thm:global-perfect} are satisfied for all step-size policies of the form $\step_{\run}\propto 1/n^{\beta}$, $\beta\in(0,1]$;
however, in the presence of errors and uncertainty, \cref{thm:global-imperfect} guarantees convergence only when $\beta\in(1/2,1]$.

The ``critical'' value $\beta = 1/2$ is tied to the finite \acl{MSE} hypothesis \eqref{eq:MSE}.
If the players' gradient observations have finite moments up to some order $q>2$, a more refined \acl{SA} argument can be used to show that \cref{thm:global-imperfect} still holds under the lighter requirement $\sum_{\run=\start}^{\infty} \step_{\run}^{1+q/2} < \infty$.
Thus, even in the presence of noise, it is possible to employ \eqref{eq:DA} with any step-size sequence of the form $\step_{\run} \propto 1/n^{\beta}$, $\beta\in(0,1]$, provided that the noise process $\noise_{\run}$ has $\exof{\dnorm{\noise_{\run+1}}^{q} \given \filter_{\run}} < \infty$ for some $q > 2/\beta - 2$.
In particular, if the noise affecting the players' observations has finite moments of all orders (for instance, if $\noise_{\run}$ is sub-exponential or sub-Gaussian), it is possible to recover essentially all the step-size policies covered by \cref{thm:global-perfect}.
\end{remark}

\subsection{Local convergence}
\label{sec:local}

The results of the previous section show that \eqref{eq:DA} converges globally to states (or sets) that are globally stable, even under noise and uncertainty.
In this section, we show that \eqref{eq:DA} remains locally convergent to states that are only locally stable with probability arbitrarily close to $1$.

For simplicity, we begin with the deterministic, perfect feedback case:

\begin{theorem}
\label{thm:local-perfect}
Suppose that \eqref{eq:DA} is run with
perfect feedback \textpar{$\noisedev=0$},
choice maps satisfying \eqref{eq:Fench-reg},
and
a sufficiently small step-size with
$\sum_{\iRun=\start}^{\run} \step_{\iRun}^{2} \big/ \sum_{\iRun=\start}^{\run} \step_{\iRun} \to 0$.
If $\eqset$ is a stable set of \aclp{NE}, there exists a neighborhood $U$ of $\eqset$ such that $\act_{\run}$ converges to $\eqset$ whenever $\act_{\start}\in U$.
\end{theorem}

\begin{Proof}
As in the proof of \cref{thm:global-perfect}, let $U_{\eps} = \setdef{x=\choice(y)}{\fench(\eqset,y) < \eps}$.
Since $\eqset$ is stable, there exists some $\eps>0$ and some $\sharp>0$ satisfying \eqref{eq:gradbound-U} and such that \eqref{eq:VS} holds throughout $U_{\eps}$.
If $\act_{\start}\in U_{\eps}$ and $\step_{\start} \leq \min\{2\sharp\strong/\vbound^{2},\sqrt{\eps\strong}/\vbound\}$, the same induction argument as in the proof of \cref{thm:global-perfect} shows that $\act_{\run}\in U_{\eps}$ for all $\run$.
Since \eqref{eq:VS} holds throughout $U_{\eps}$, \cref{lem:visit} shows that $\act_{\run}$ visits any neighborhood of $\eqset$ infinitely many times.
Thus, by the same argument as in the proof of \cref{thm:global-perfect}, we get $\act_{\run}\to\eqset$.
\end{Proof}

The key idea in the proof of \cref{thm:local-perfect} is that if the step-size of \eqref{eq:DA} is small enough, $\act_{\run} = \choice(\score_{\run})$ always remains within the ``basin of attraction'' of $\eqset$;
hence, local convergence can be obtained in the same way as global convergence for a game with smaller action spaces.
However, if the players' feedback is subject to estimation errors and uncertainty, a single unlucky instance could drive $\act_{\run}$ away from said basin, possibly never to return.
Consequently, any local convergence result in the presence of noise is necessarily probabilistic in nature.

Conditioning on the event that $\act_{\run}$ stays close to $\eqset$, local convergence can be obtained as in the proof of \cref{thm:global-imperfect}.
Nevertheless, showing that this event occurs with controllably high probability requires a completely different analysis.
This is the essence of our next result:

\begin{theorem}
\label{thm:local-imperfect}
Fix a confidence level $\delta>0$ and suppose that \eqref{eq:DA} is run with a sufficiently small step-size $\step_{\run}$ satisfying $\sum_{\run=\start}^{\infty} \step_{\run}^{2} < \infty$ and $\sum_{\run=\start}^{\infty} \step_{\run} = \infty$.
If $\eqset$ is stable and \crefrange{eq:zeromean}{eq:Lipschitz} hold, then $\eqset$ is locally attracting with probability at least $1-\delta$;
more precisely, there exists a neighborhood $U$ of $\eqset$ such that
\begin{equation}
\probof{\act_{\run}\to\eqset \given \act_{\start}\in U}
	\geq 1-\delta.
\end{equation}
\end{theorem}

\begin{corollary}
\label{cor:local}
Let $\eq$ be a \acl{NE} with negative-definite Hessian matrix $\hessmat(\eq) \ml 0$.
Then, with assumptions as above, $\eq$ is locally attracting with probability arbitrarily close to $1$.
\end{corollary}

\begin{Proof}[Proof of \cref{thm:local-imperfect}]
Let $U_{\eps} = \setdef{x=\choice(y)}{\fench(\eqset,y) < \eps}$ and pick $\eps>0$ small enough so that \eqref{eq:VS} holds for all $x\in U_{3\eps}$.
Assume further that $\act_{\start}\in U_{\eps}$ so there exists some $\eq\in\eqset$ such that $\fench(\eq,\score_{\start}) < \eps$.
Then,
for all $\run$, \cref{prop:Fenchel} yields
\begin{flalign}
\label{eq:Dbound-stoch}
\fench(\eq,\score_{\run+1})
	&\leq
	\fench(\eq,\score_{\run})
	+ \step_{\run} \braket{\payv(\act_{\run})}{\act_{\run} - \eq}
	+ \step_{\run} \snoise_{\run+1}
	+ \frac{\step_{\run}^{2}}{2\strong} \dnorm{\est_{\run+1}}^{2},
\end{flalign}
where we have set $\snoise_{\run+1} = \braket{\noise_{\run+1}}{\act_{\run} - \eq}$.

We first claim that $\sup_{\run} \sum_{\iRun=\start}^{\run} \step_{\iRun} \snoise_{\iRun+1} \leq \eps$ with probability at least $1-\delta/2$ if $\step_{\run}$ is chosen appropriately.
Indeed, set $S_{\run+1} = \sum_{\iRun=\start}^{\run} \step_{\iRun} \snoise_{\iRun+1}$ and let $E_{\run,\eps}$ denote the event $\{\sup_{1\leq \iRun\leq \run+1} \abs{S_{\iRun}} \geq \eps\}$.
Since $S_{\run}$ is a martingale, Doob's maximal inequality \citep[Theorem 2.1]{HH80} yields
\begin{equation}
\label{eq:prob-bound1}
\probof{E_{\run+1,\eps}}
	\leq \frac{\exof{\abs{S_{\run+1}}^{2}}}{\eps^{2}}
	\leq \frac{\noisevar \norm{\feas}^{2} \sum_{\iRun=\start}^{\run} \step_{\iRun}^{2}}{\eps^{2}},
\end{equation}
where we used the variance estimate
\begin{flalign}
\exof{\snoise_{\iRun+1}^{2}}
	&= \exof{\exof{\abs{\braket{\noise_{\iRun+1}}{\act_{\iRun} - \eq}}^{2} \given \filter_{\iRun}}}
	\notag\\
	&\leq \exof{\exof{\dnorm{\noise_{\iRun+1}}^{2} \norm{\act_{\iRun} - \eq}^{2} \given \filter_{\iRun}}}
	\leq \noisevar \norm{\feas}^{2},
\end{flalign}
and the fact that $\exof{\snoise_{\iRun+1}\snoise_{\ell+1}} = \exof{\exof{\snoise_{\iRun+1} \snoise_{\ell+1}} \given \filter_{\iRun\vee\ell}} = 0$ whenever $\iRun\neq\ell$.
Since $E_{\run+1,\eps} \supseteq E_{\run,\eps}\supseteq\dotsc$, the event $E_{\eps} = \union_{\run=\start}^{\infty} E_{\run,\eps}$ occurs with probability $\prob(E_{\eps}) \leq \Gamma_{2} \noisevar \norm{\feas}^{2} / \eps^{2}$, where $\Gamma_{2} \equiv \sum_{\run=\start}^{\infty} \step_{\run}^{2}$.
Thus, if $\Gamma_{2} \leq \delta\eps^{2}/(2\noisevar\norm{\feas}^{2})$, we get $\probof{E_{\eps}} \leq \delta/2$.

We now claim that the process $R_{\run+1} = (2\strong)^{-1} \sum_{\iRun=\start}^{\run} \step_{\iRun}^{2} \dnorm{\est_{\iRun+1}}^{2}$ is also bounded from above by $\eps$ with probability at least $1-\delta/2$ if $\step_{\run}$ is chosen appropriately.
Indeed, working as above, let $F_{\run,\eps}$ denote the event $\{\sup_{1\leq \iRun\leq \run+1} R_{\iRun} \geq \eps\}$.
Since $R_{\run}$ is a nonnegative submartingale, Doob's maximal inequality again yields
\begin{equation}
\label{eq:prob-bound2}
\probof{F_{\run+1,\eps}}
	\leq \frac{\exof{R_{\run+1}}}{\eps}
	\leq \frac{\vbound^{2} \sum_{\iRun=\start}^{\run} \step_{\iRun}^{2}}{2\strong\eps}.
\end{equation}
Consequently, the event $F_{\eps} = \union_{\run=\start}^{\infty} F_{\run,\eps}$ occurs with probability $\probof{F_{\eps}} \leq \Gamma_{2} \vbound^{2}/\eps \leq \delta/2$ if $\step_{\run}$ is chosen so that $\Gamma_{2} \leq \strong\delta\eps/\vbound^{2}$.

Assume therefore that $\Gamma_{2} \leq \min\{\delta\eps^{2}/(2\noisevar\norm{\feas}^{2}),\,\strong\delta\eps/\vbound^{2}\}$.
The above shows that $\probof{\bar E_{\eps} \cap \bar F_{\eps}} = 1 - \probof{E_{\eps}\cup F_{\eps}} \geq 1 - \delta/2 - \delta/2 = 1 - \delta$, i.e. $S_{\run}$ and $R_{\run}$ are both bounded from above by $\eps$ for all $\run$ and all $\eq$ with probability at least $1-\delta$.
Since $\fench(\eq,\score_{\start}) \leq \eps$ by assumption, we readily get $\fench(\eq,\score_{\start}) \leq 3\eps$ if $\bar E_{\eps}$ and $\bar F_{\eps}$ both hold.
Furthermore, telescoping \eqref{eq:Dbound-stoch} yields
\begin{equation}
\label{eq:Dbound-n}
\fench(\eq,\score_{\run+1})
	\leq \fench(\eq,\score_{\start})
	+ \sum_{\iRun=\start}^{\run} \braket{\payv(\act_{\iRun})}{\act_{\iRun} - \eq}
	+ S_{\run+1}
	+ R_{\run+1}
	\quad
	\text{for all $\run$},
\end{equation}
so if we assume inductively that $\fench(\eq,\score_{\iRun}) \leq 3\eps$ for all $\iRun\leq\run$ (implying that $\braket{\payv(\act_{\iRun})}{\act_{\iRun} - \eq} \leq 0$ for all $k\leq n$), we also get $\fench(\eq,\score_{\run+1}) \leq 3\eps$ if neither $E_{\eps}$ nor $F_{\eps}$ occur.
Since $\probof{E_{\eps}\cup F_{\eps}} \leq \delta$, we conclude that $\act_{\run}$ stays in $U_{3\eps}$ for all $\run$ with probability at least $1-\delta$.
In turn, when this is the case, \cref{lem:visit} shows that $\eqset$ is recurrent under $\act_{\run}$.
Hence, by repeating the same steps as in the proof of \cref{thm:global-imperfect}, we get $\act_{\run}\to\eqset$ with probability at least $1-\delta$, as claimed.
\end{Proof}

\subsection{Convergence in zero-sum concave games}
\label{sec:zerosum}

We close this section by examining the asymptotic behavior of \eqref{eq:DA} in $2$-player, concave-convex zero-sum games.
To do so, let $\players = \{\playOne,\playTwo\}$ denote the set of players with corresponding payoff functions $\pay_{\playOne} = - \pay_{\playTwo}$ respectively concave in $x_{\playOne}$ and $x_{\playTwo}$.
Letting $\pay \equiv \pay_{\playOne} = -\pay_{\playTwo}$, the \emph{value} of the game is defined as
\begin{equation}
\label{eq:value}
\val
	= \max_{x_{\playOne}\in\feas_{\playOne}} \min_{x_{\playTwo}\in\feas_{\playTwo}} \pay(x_{\playOne},x_{\playTwo})
	= \min_{x_{\playTwo}\in\feas_{\playTwo}} \max_{x_{\playOne}\in\feas_{\playOne}} \pay(x_{\playOne},x_{\playTwo}).
\end{equation}
The solutions of the concave-convex saddle-point problem \eqref{eq:value} are the \aclp{NE} of $\game$ and the players' equilibrium payoffs are $\val$ and $-\val$ respectively.

In the ``perfect feedback'' case ($\noisedev = 0$), \cite{Nes09} showed that the ergodic average
\begin{equation}
\label{eq:ergodic}
\bar \act_{\run}
	= \frac{\sum_{\iRun=\start}^{\run} \step_{\iRun} \act_{\iRun}}{\sum_{\iRun=\start}^{\run} \step_{\iRun}}
\end{equation}
of the sequence of play generated by \eqref{eq:DA} converges to equilibrium.
With imperfect feedback and steep $h$,%
\footnote{When $h$ is steep, the \acl{MD} algorithm examined by \cite{NJLS09} is a special case of the \acl{DA} method of \cite{Nes09}.
This is no longer the case if $h$ is not steep, so the analysis of \cite{NJLS09} does not apply to \eqref{eq:DA}.
In the online learning literature, this difference is sometimes referred to as ``greedy'' vs. ``lazy'' \acl{MD}.}
\cite{NJLS09} further showed that $\bar \act_{\run}$ converges in expectation to the game's set of \aclp{NE}, provided that \eqref{eq:zeromean} and \eqref{eq:MSE} hold.
Our next result provides an almost sure version of this result which is also valid for nonsteep $h$:

\begin{theorem}
\label{thm:zerosum}
Let $\game$ be a concave $2$-player zero-sum game.
If \eqref{eq:DA} is run with
imperfect feedback satisfying \crefrange{eq:zeromean}{eq:MSE}
and
a step-size $\step_{\run}$ such that
$\sum_{\run=\start}^{\infty} \step_{\run}^{2} < \infty$ and $\sum_{\run=\start}^{\infty} \step_{\run} = \infty$,
the ergodic average $\bar \act_{\run}$ of $\act_{\run}$ converges to the set of \aclp{NE} of $\game$ \as.
\end{theorem}


\begin{Proof}[Proof of \cref{thm:zerosum}]
Consider the gap function
\begin{equation}
\gap(x)
	= \val - \min_{\base_{\playTwo}\in\feas_{\playTwo}} \pay(x_{\playOne},\base_{\playTwo})
	+ \max_{\base_{\playOne}\in\feas_{\playOne}} \pay(\base_{\playOne},x_{\playTwo}) - \val
	= \max_{\base\in\feas} \sum_{\play\in\players} \pay_{\play}(\base_{\play};x_{-\play}).
\end{equation}
Obviously, $\gap(x) \geq 0$ with equality if and only if $x$ is a \acl{NE}, so it suffices to show that $\gap(\bar \act_{\run})\to 0$ \as.

To do so, pick some $\base\in\feas$.
Then, as in the proof of \cref{thm:global-imperfect}, we have
\begin{equation}
\label{eq:Dn}
\fench(\base,\score_{\run+1})
	\leq \fench(\base,\score_{\run})
	+ \step_{\run} \braket{\payv(\act_{\run})}{\act_{\run} - \base}
	+ \step_{\run} \snoise_{\run+1}
	+ \frac{1}{2\strong} \step_{\run}^{2} \dnorm{\est_{\run+1}}^{2}.
\end{equation}
Hence, after rearranging and telescoping, we get
\begin{equation}
\label{eq:aux1}
\sum_{\iRun=\start}^{\run} \step_{\iRun} \braket{\payv(\act_{\iRun})}{\base - \act_{\iRun}}
	\leq \fench(\base,\score_{\start})
	+ \sum_{\iRun=\start}^{\run} \step_{\iRun} \snoise_{\iRun+1}
	+ \frac{1}{2\strong} \sum_{\iRun=\start}^{\run} \step_{\iRun}^{2} \dnorm{\est_{\iRun+1}}^{2},
\end{equation}
where $\snoise_{\run+1} = \braket{\noise_{\run+1}}{\act_{\run} - \base}$ and we used the fact that $\fench(\base,\score_{\run})\geq0$.
By concavity, we also have
\begin{equation}
\label{eq:aux2}
\braket{\payv(x)}{\base - x}
	= \sum_{\play\in\players} \braket{\payv_{\play}(x)}{\base_{\play} - x_{\play}}
	\geq \sum_{\play\in\players} \left[ \pay_{\play}(\base_{\play};x_{-\play}) - \pay_{\play}(x) \right]
	= \sum_{\play\in\players} \pay_{\play}(\base_{\play};x_{-\play}),
\end{equation}
for all $x\in\feas$.
Therefore, letting $\tau_{\run} = \sum_{\iRun=\start}^{\run} \step_{\iRun}$, we get
\begin{flalign}
\label{eq:aux3}
\frac{1}{\tau_{\run}} \sum_{\iRun=\start}^{\run} \step_{\iRun} \braket{\payv(\act_{\iRun})}{\base - \act_{\iRun}}
	&\geq \frac{1}{\tau_{\run}} \sum_{\iRun=\start}^{\run} \step_{\iRun} \sum_{\play\in\players} \pay_{\play}(\base_{\play};\act_{-\play,\iRun})
	\notag\\
	&\geq \pay(\base_{\playOne},\bar \act_{\playTwo,\run}) - \pay(\bar \act_{\playOne,\run},\base_{\playTwo})
	\notag\\
	&= \sum_{\play\in\players} \pay_{\play}(\base_{\play};\bar \act_{-\play,\run}),
\end{flalign}
where we used the fact that $\pay$ is concave-convex in the second line.
Thus, combining \eqref{eq:aux1} and \eqref{eq:aux3}, we finally obtain
\begin{equation}
\label{eq:aux4}
\sum_{\play\in\players} \pay_{\play}(\base_{\play};\bar \act_{-\play,\run})
	\leq \frac{\fench(\base,\score_{\start}) + \sum_{\iRun=\start}^{\run} \step_{\iRun}\snoise_{\iRun+1} + (2\strong)^{-1} \sum_{\iRun=\start}^{\run} \step_{\iRun}^{2} \dnorm{\est_{\iRun+1}}^{2}}{\tau_{\run}}.
\end{equation}

As before, the \acl{LLN} \citep[Theorem 2.18]{HH80} yields $\tau_{\run}^{-1} \sum_{\iRun=\start}^{\run} \step_{\iRun} \snoise_{\iRun+1} \to 0$ \as.
Furthermore, given that $\exof{\dnorm{\est_{\run+1}}^{2} \given \filter_{\run}} \leq \vbound^{2}$ and $\sum_{\iRun=\start}^{\run} \step_{\iRun}^{2} < \infty$, we also get $\tau_{\run}^{-1} \sum_{\iRun=\start}^{\run} \step_{\iRun}^{2} \dnorm{\est_{\iRun+1}}^{2} \to 0$ by Doob's martingale convergence theorem \citep[Theorem~2.5]{HH80},
implying in turn that $\sum_{\play\in\players} \pay_{\play}(\base_{\play};\bar \act_{-\play,\run}) \to 0$ \as.
Since $\base$ is arbitrary, we conclude that $\gap(\bar \act_{\run}) \to 0$ \as, as claimed.
\end{Proof}

%% file: Tables/Hypotheses.tex

\begin{tabular}{cll}
	&\sc{Hypothesis}
	&\sc{Precise statement}
	\\
\hline
\eqref{eq:zeromean}
	&Zero-mean errors
	&$\exof{\noise_{\run+1} \given \filter_{\run}} = 0$
	\\
\hline
\eqref{eq:MSE}
	&Finite error variance
	&$\exof{\dnorm{\noise_{\run+1}}^{2} \given \filter_{\run}} \leq \noisevar$
	\\
\hline
\eqref{eq:Fench-reg}
	&Bregman reciprocity
	&$\fench(\base,y_{\run})\to0$ whenever $\choice(y_{\run})\to\base$
	\\
\hline
\eqref{eq:Lipschitz}
	&Lipschitz gradients
	&$\payv(x)$ is Lipschitz continuous
	\\
\hline
\end{tabular}

%% file: Finite.tex

As a concrete application of the analysis of the previous section, we turn to the asymptotic behavior of \eqref{eq:DA} in \emph{finite} games.
Briefly recalling the setup of \cref{ex:finite}, each player in a finite game $\fingame\equiv\fingamefull$ chooses a pure strategy $\pure_{\play}$ from a finite set $\pures_{\play}$ and receives a payoff of $\pay_{\play}(\pure_{1},\dotsc,\pure_{\nPlayers})$.
Pure strategies are drawn based on the players' mixed strategies $x_{\play}\in\feas_{\play}\equiv\simplex(\pures_{\play})$, so each player's expected payoff is given by the multilinear expression \eqref{eq:pay-mixed}.
Accordingly, the individual payoff gradient of player $\play\in\players$ in the mixed profile $x = (x_{1},\dotsc,x_{\nPlayers})$ is the (mixed) payoff vector
\(
\payv_{\play}(x)
	= \grad_{x_{\play}} \pay_{\play}(x_{\play};x_{-\play})
	= (\pay_{\play}(\pure_{\play};x_{-\play}))_{\pure_{\play}\in\pures_{\play}}
\)
of \cref{eq:payv-finite}.

Consider now the following learning scheme:
At stage $\run$, every player $\play\in\players$ selects a pure strategy $\pure_{\play,\run}\in\pures_{\play}$ according to their individual mixed strategy $\act_{\play,\run}\in\feas_{\play}$.
Subsequently, each player observes \textendash\ or otherwise calculates \textendash\ the payoffs of their pure strategies $\pure_{\play}\in\pures_{\play}$ against the chosen actions $\pure_{-\play,\run}$ of all other players (possibly subject to some random estimation error).
Specifically, we posit that each player receives as feedback the ``noisy'' payoff vector
\begin{equation}
\label{eq:payv-finite-noise}
\est_{\play,\run+1}
	= (\pay_{\play}(\pure_{\play};\pure_{-\play,\run}))_{\pure_{\play}\in\pures_{\play}}
	+ \noise_{\play,\run+1},
\end{equation}
where the error process $\noise_{\run} = (\noise_{\play,\run})_{\play\in\players}$ is assumed to satisfy \labelcref{eq:zeromean,eq:MSE}.
Then, based on this feedback, players update their mixed strategies and the process repeats (for a concrete example, see \cref{alg:logit}).

\smallskip

In the rest of this section, we study the long-term behavior of this adaptive learning process.
Specifically, we focus on:
\begin{inparaenum}
[\itshape a\upshape)]
\item
the elimination of dominated strategies;
\item
convergence to strict \aclp{NE};
and
\item
convergence to equilibrium in $2$-player, zero-sum games.
\end{inparaenum}

\input{Logit}

\subsection{Dominated strategies}
\label{sec:dominated}

We say that a pure strategy $\pure_{\play}\in\pures_{\play}$ of a finite game $\fingame$ is \emph{dominated} by $\purealt_{\play}\in\pures_{\play}$ (and we write $\pure_{\play}\ml\purealt_{\play}$) if
\begin{equation}
\label{eq:dominated}
\txs
\pay_{\play}(\pure_{\play};x_{-\play})
	< \pay_{\play}(\purealt_{\play};x_{-\play})
	\quad
	\text{for all $x_{-\play}\in\feas_{-\play}\equiv\prod_{\playalt\neq\play}\feas_{\playalt}$}.
\end{equation}
Put differently, $\pure_{\play} \ml \purealt_{\play}$ if and only if $\payv_{\play\pure_{\play}}(x) < \payv_{\play\purealt_{\play}}(x)$ for all $x\in\feas$.
In turn, this implies that the payoff gradient of player $\play$ points consistently towards the face $x_{\play\pure_{\play}} = 0$ of $\feas_{\play}$, so it is natural to expect that $\pure_{\play}$ is eliminated under \eqref{eq:DA}.
Indeed, we have:

\begin{theorem}
\label{thm:dominated}
Suppose that \eqref{eq:DA} is run with
noisy payoff observations of the form \eqref{eq:payv-finite-noise}
and
a step-size sequence $\step_{\run}$ satisfying \eqref{eq:step-loose}.
If $\pure_{\play}\in\pures_{\play}$ is dominated, then $\act_{\play\pure_{\play},\run}\to0$ \as.
\end{theorem}

\begin{Proof}
Suppose that $\pure_{\play} \ml \purealt_{\play}$ for some $\purealt_{\play}\in\pures_{\play}$.
Then, suppressing the player index $\play$ for simplicity, \eqref{eq:DA} gives
\begin{flalign}
\label{eq:ydiff}
\score_{\purealt,\run+1} - \score_{\pure,\run+1}
	&= \sharp_{\purealt\pure}
	+ \sum_{\iRun=\start}^{\run} \step_{\iRun} \left[ \est_{\purealt,\iRun+1} - \est_{\pure,\iRun+1} \right]
	\notag\\
	&= \sharp_{\purealt\pure}
	+ \sum_{\iRun=\start}^{\run} \step_{\iRun} \left[ \payv_{\purealt}(\act_{\iRun}) - \payv_{\pure}(\act_{\iRun}) \right]
	+ \sum_{\iRun=\start}^{\run} \step_{\iRun} \zeta_{\iRun+1},
\end{flalign}
where we set
$\sharp_{\purealt\pure} = \score_{\purealt,\start} - \score_{\pure,\start}$
and
\begin{equation}
\zeta_{\iRun+1}
	= \exof{\est_{\purealt,\iRun+1} - \est_{\pure,\iRun+1} \given \filter_{\iRun}} - [\payv_{\purealt}(\act_{\iRun}) - \payv_{\pure}(\act_{\iRun})].
\end{equation}
Since $\pure\ml\purealt$, there exists some $\sharp>0$ such that $\payv_{\purealt}(x) - \payv_{\pure}(x) \geq \sharp$ for all $x\in\feas$.
Then, \eqref{eq:ydiff} yields
\begin{equation}
\score_{\purealt,\run+1} - \score_{\pure,\run+1}
	\geq \sharp_{\purealt\pure} + \tau_{\run} \left[ \sharp + \frac{\sum_{\iRun=\start}^{\run} \step_{\iRun} \zeta_{\iRun+1}}{\tau_{\run}} \right],
\end{equation}
where $\tau_{\run} = \sum_{\iRun=\start}^{\run} \step_{\iRun}$.
As in the proof of \cref{thm:terminal}, the \acl{LLN} for \aclp{MDS} \citefp[Theorem 2.18]{HH80} implies that $\tau_{\run}^{-1}\sum_{\iRun=\start}^{\run} \step_{\iRun} \zeta_{\iRun+1} \to 0$ under the step-size assumption \eqref{eq:step-loose}, so $\score_{\purealt,\run} - \score_{\pure,\run} \to \infty$ \as.

Suppose now that $\limsup_{\run\to\infty} \act_{\pure,\run} = 2\eps$ for some $\eps>0$.
By descending to a subsequence if necessary, we may assume that $\act_{\pure,\run} \geq \eps$ for all $\run$, so if we let $\act_{\run}' = \act_{\run} + \eps (\bvec_{\purealt} - \bvec_{\pure})$, the definition of $\choice$ gives
\begin{equation}
h(\act_{\run}')
	\geq h(\act_{\run}) + \braket{\score_{\run}}{\act_{\run}' - \act_{\run}}
	= h(\act_{\run}) + \eps \parens{\score_{\purealt,\run} - \score_{\pure,\run}}
	\to \infty,
\end{equation}
a contradiction.
This implies that $\act_{\pure,\run}\to0$ \as, as asserted.
\end{Proof}

\subsection{Strict equilibria}
\label{sec:strict}

A \acl{NE} $\eq$ of a finite game is called \emph{strict} when \eqref{eq:Nash} holds as a strict inequality for all $x_{\play}\neq\eq_{\play}$, i.e. when no player can deviate unilaterally from $\eq$ without \emph{reducing} their payoff (or, equivalently, when every player has a unique best response to $\eq$).
This implies that strict \aclp{NE} are pure strategy profiles $\eq=(\peq_{1},\dotsc,\peq_{\nPlayers})$ such that
\begin{equation}
\label{eq:strict-finite}
\pay_{\play}(\peq_{\play};\peq_{-\play})
	> \pay_{\play}(\pure_{\play};\peq_{-\play})
	\quad
	\text{for all $\pure_{\play}\in\pures_{\play}\setminus\{\peq_{\play}\}$, $\play\in\players$}.
\end{equation}
Strict \aclp{NE} can be characterized further as follows:

\begin{proposition}
\label{prop:strict}
Then, the following are equivalent:
\begin{enumerate}
[\indent a\upshape)]
\item
$\eq$ is a strict \acl{NE}.
\item
$\braket{\payv(\eq)}{z} \leq 0$ for all $z\in\tcone(\eq)$ with equality if and only if $z=0$.
\item
$\eq$ is stable.
\end{enumerate}
\end{proposition}

Thanks to the above characterization of strict equilibria (proven in \cref{app:aux}), the convergence analysis of \cref{sec:analysis} yields:

\begin{proposition}
\label{prop:conv-strict}
Let $\eq$ be a strict equilibrium of a finite game $\fingame$.
Suppose further that \eqref{eq:DA} is run with noisy payoff observations of the form \eqref{eq:payv-finite-noise} and a sufficiently small step-size $\step_{\run}$ such that
$\sum_{\run=\start}^{\infty} \step_{\run}^{2} < \infty$ and $\sum_{\run=\start}^{\infty} \step_{\run} = \infty$.
If \crefrange{eq:zeromean}{eq:Fench-reg} hold, $\eq$ is locally attracting with arbitrarily high probability;
specifically, for all $\delta>0$, there exists a neighborhood $U$ of $\eq$ such that
\begin{equation}
\probof{\act_{\run}\to\eq \given \act_{\start}\in U}
	\geq 1-\delta.
\end{equation}
\end{proposition}

\begin{Proof}
We first show that $\exof{\est_{\run+1} \given \filter_{\run}} = \payv(\act_{\run})$.
Indeed, for all $\play\in\players$, $\pure_{\play}\in\pures_{\play}$, we have
\begin{equation}
\exof{\est_{\play\pure_{\play},\run+1} \given \filter_{\run}}
	= \sum_{\pure_{-\play}\in\pures_{-\play}}\!\!\!
	\pay_{\play}(\pure_{\play};\pure_{-\play})\,
	\act_{\pure_{-\play},\run}
	+ \exof{\noise_{\play\pure_{\play},\run+1} \given \filter_{\run}}
	= \pay_{\play}(\pure_{\play};\act_{-\play,\run}),
\end{equation}
where, in a slight abuse of notation, we set $\act_{\pure_{-\play},\run}$ for the joint probability assigned to the pure strategy profile $\pure_{-\play}$ of all players other than $\play$ at stage $\run$.

By \eqref{eq:payv-finite}, it follows that $\exof{\est_{\run+1} \given \filter_{\run}} = \payv(\act_{\run})$ so the estimator \eqref{eq:payv-finite-noise} is unbiased in the sense of \labelcref{eq:zeromean}.
Hypothesis \eqref{eq:MSE} can be verified similarly, so the estimator \eqref{eq:payv-finite-noise} satisfies \eqref{eq:estimates}.
Since $\eq$ is stable by \cref{prop:strict} and $\payv(x)$ is multilinear (so \eqref{eq:Lipschitz} is satisfied automatically), our assertion follows from \cref{thm:local-imperfect}.
\end{Proof}

In the special case of logit-based learning (\cref{ex:logit}), \citef{CHM17-SAGT} showed that \cref{alg:logit}
converges locally to strict \aclp{NE} under similar information assumptions.
\cref{prop:strict} essentially extends this result to the entire class of regularized learning processes induced by \eqref{eq:DA} in finite games, showing that the logit choice map \eqref{eq:choice-logit} has no special properties in this regard.
\citef{CHM17-SAGT} further showed that the convergence rate of logit-based learning is exponential in the algorithm's ``running horizon'' $\tau_{\run} = \sum_{\iRun=\start}^{\run} \step_{\iRun}$.
This rate is closely linked to the logit choice model, and different choice maps yield different convergence speeds;
we discuss this issue in more detail in \cref{sec:rates}.

\subsection{Convergence in zero-sum games}
\label{sec:zerosum-fin}

We close this section with a brief discussion of the ergodic convergence properties of \eqref{eq:DA} in finite two-player zero-sum games.
In this case, the analysis of \cref{sec:zerosum} readily yields:

\begin{corollary}
\label{cor:zerosum}
Let $\fingame$ be a finite $2$-player zero-sum game.
If \eqref{eq:DA} is run with
noisy payoff observations of the form \eqref{eq:payv-finite-noise}
and
a step-size $\step_{\run}$ such that
$\sum_{\run=\start}^{\infty} \step_{\run}^{2} < \infty$ and $\sum_{\run=\start}^{\infty} \step_{\run} = \infty$,
the ergodic average $\bar \act_{\run} = \sum_{\iRun=\start}^{\run} \step_{\iRun} \act_{\iRun} \big/ \sum_{\iRun=\start}^{\run} \step_{\iRun}$ of the players' mixed strategies converges to the set of \aclp{NE} of $\fingame$ \as.
\end{corollary}

\begin{Proof}
As in the proof of \cref{prop:conv-strict}, the estimator \eqref{eq:payv-finite-noise} satisfies $\exof{\est_{\run+1} \given \filter_{\run}} = \payv(\act_{\run})$, so \labelcref{eq:zeromean,eq:MSE} also hold in the sense of \eqref{eq:estimates}.
Our claim then follows from \cref{thm:zerosum}.
\end{Proof}

\begin{remark}
In a very recent paper, \citef{BM17} showed that the time average $\bar\act(t) = t^{-1}\int_{0}^{t} \act(s) \dd s$ of the players' mixed strategies under \eqref{eq:DA-cont} with Brownian payoff shocks converges to \acl{NE} in $2$-player, zero-sum games.
\cref{cor:zerosum} may be seen as a discrete-time version of this result.
\end{remark}

%% file: Logit.tex

\begin{algorithm}[tbp]
\caption{Logit-based learning in finite games (\cref{ex:logit}).}
\label{alg:logit}
\normalsize
\begin{algorithmic}[1]
	\REQUIRE
	step-size sequence $\step_{\run}\propto 1/\run^{\beta}$, $\beta\in(0,1]$;
	initial scores $\score_{\play}\in\R^{\pures_{\play}}$
	\FOR{$\run=\start,2,\dotsc$}
		\FOR{every player $\play\in\players$}
			\STATE set $\act_{\play} \leftarrow \logit_{\play}(\score_{\play})$;
			\hfill
			\COMMENT{mixed strategy}
			\\[2pt]
			\STATE play $\pure_{\play} \sim \act_{\play}$;
			\hfill
			\COMMENT{choose action}
			\\[2pt]
			\STATE observe $\est_{\play}$;
			\hfill
			\COMMENT{estimate payoffs}
			\\[2pt]
			\STATE update $\score_{\play} \leftarrow \score_{\play} + \step_{\run} \est_{\play}$;
			\hfill
			\COMMENT{update scores}
		\ENDFOR
	\ENDFOR
\end{algorithmic}
\end{algorithm}

%% file: Rates.tex

\subsection{Ergodic convergence rate}
\label{sec:rate-ergodic}

In this section, we focus on the rate of convergence of \eqref{eq:DA} to stable equilibrium states (and/or sets thereof).
To that end, we will measure the speed of convergence to a globally stable set $\eqset\subseteq\feas$ via the \emph{equilibrium gap function}
\begin{equation}
\label{eq:gap}
\gap(x)
	= \inf_{\eq\in\eqset} \braket{\payv(x)}{\eq - x}.
\end{equation}
By \cref{def:VS}, $\gap(x) \geq 0$ with equality if and only if $x\in\eqset$, so $\gap(x)$ can be seen as a (game-dependent) measure of the distance between $x$ and the target set $\eqset$.
This can be seen more clearly in the case of \emph{strongly} stable equilibria, defined here as follows:

\begin{definition}
\label{def:SVS}
We say that $\eq\in\feas$ is \emph{strongly stable} if there exists some $L>0$ such that
\begin{alignat}{3}
\label{eq:SVS}
\braket{\payv(x)}{x - \eq}
	&\leq - L \norm{x - \eq}^{2}
	&\quad
	&\text{for all $x\in\feas$}.
	\\
\intertext{More generally, a closed subset $\eqset$ of $\feas$ is called strongly stable if}
\label{eq:SVS-set}
\txs
\braket{\payv(x)}{x - \eq}
	&\leq - L \dist(\eqset,x)^{2}
	&\quad
	&\text{for all $x\in\feas$, $\eq\in\eqset$}.
\end{alignat}
\end{definition}

Obviously, $\gap(x) \geq L \dist(\eqset,x)^{2}$ if $\eqset$ is $L$-strongly stable, i.e. $\gap(x)$ grows at least quadratically near strongly stable sets \textendash\ just like strongly convex functions grow quadratically around their minimum points.
With this in mind, we provide below an explicit estimate for the decay rate of the average equilibrium gap $\bar\gap_{\run} = \sum_{\iRun=\start}^{\run} \step_{\iRun}\gap(\act_{\iRun})\big/\sum_{\iRun=\start}^{\run} \step_{\iRun}$ in the spirit of \cite{NJLS09}:

\begin{theorem}
\label{thm:gap-ergodic}
Suppose that \eqref{eq:DA} is run with imperfect gradient information satisfying \crefrange{eq:zeromean}{eq:MSE}.
Then
\begin{equation}
\label{eq:gap-ergodic}
\exof{\bar\gap_{\run}}
	\leq \frac{\fench_{\start} + \vbound^{2}/(2\strong) \sum_{\iRun=\start}^{\run} \step_{\iRun}^{2}}{\sum_{\iRun=\start}^{\run} \step_{\iRun}},
\end{equation}
where $\fench_{\start} = \fench(\eqset,\score_{\start})$.
If, in addition, $\sum_{\run=\start}^{\infty} \step_{\run}^{2} < \infty$,
we have
\begin{equation}
\label{eq:gap-ergodic-as}
\bar\gap_{\run}
	\leq \frac{A}{\sum_{\iRun=\start}^{\run} \step_{\iRun}}
	\quad
	\text{for all $\run$ \as},
\end{equation}
where $A>0$ is a finite random variable such that, with probability at least $1-\delta$,
\begin{equation}
\label{eq:gap-constant}
A
	\leq \fench_{\start}
	+ \noisedev\norm{\feas}\kappa
	+ \kappa^{2} \vbound^{2},
\end{equation}
where $\kappa^{2} = 2 \delta^{-1} \sum_{\run=\start}^{\infty} \step_{\run}^{2}$.
\end{theorem}


\begin{corollary}
\label{cor:gap-ergodic-opt}
Suppose that \eqref{eq:DA} is initialized at $\score_{\start}=0$ and is run for $\run$ iterations with constant step-size $\step = \vbound^{-1} \sqrt{2\strong\depth/\run}$ where $\depth = \max h - \min h$.
Then,
\begin{equation}
\label{eq:gap-ergodic-opt}
\exof{\bar\gap_{\run}}
	\leq 2\vbound \sqrt{\depth/(\strong\run)}.
\end{equation}
In addition, if $\eqset$ is $L$-strongly stable,
the long-run average distance to equilibrium $\bar r_{\run} = \sum_{\iRun=\start}^{\run} \dist(\eqset,\act_{\run}) \big/ \sum_{\iRun=\start}^{\run} \step_{\iRun}$ satisfies
\begin{equation}
\exof{\bar r_{\run}}
	\leq \sqrt[4]{4 L^{-2} \vbound^{2} \depth/(\strong\run)}.
\end{equation}
\end{corollary}


\begin{Proof}[Proof of \cref{thm:gap-ergodic}]
Let $\eq\in\eqset$.
Rearranging \eqref{eq:Dbound-stoch} and telescoping yields
\begin{equation}
\label{eq:gap-ergodic2}
\sum_{\iRun=\start}^{\run} \step_{\iRun} \braket{\payv(\act_{\iRun})}{\eq - \act_{\iRun}}
	\leq \fench(\eq,\score_{\start})
	+ \sum_{\iRun=\start}^{\run} \step_{\iRun} \snoise_{\iRun+1}
	+ \frac{1}{2\strong} \sum_{\iRun=\start}^{\run} \step_{\iRun}^{2} \dnorm{\est_{\iRun+1}}^{2},
\end{equation}
where $\snoise_{\iRun+1} = \braket{\noise_{\iRun+1}}{\act_{\iRun} - \eq}$.
Thus, taking expectations on both sides, we obtain
\begin{equation}
\label{eq:gap-ergodic3}
\sum_{\iRun=\start}^{\run} \step_{\iRun} \exof{\braket{\payv(\act_{\iRun})}{\eq - \act_{\iRun}}}
	\leq \fench(\eq,\score_{\start})
	+ \frac{\vbound^{2}}{2\strong} \sum_{\iRun=\start}^{\run} \step_{\iRun}^{2}.
\end{equation}
Subsequently, minimizing both sides of \eqref{eq:gap-ergodic3} over $\eq\in\eqset$ yields
\begin{equation}
\label{eq:gap-ergodic4}
\sum_{\iRun=\start}^{\run} \step_{\iRun} \exof{\gap(\act_{\iRun})}
	\leq \fench_{\start}
	+ \frac{\vbound^{2}}{2\strong} \sum_{\iRun=\start}^{\run} \step_{\iRun}^{2},
\end{equation}
where we used Jensen's inequality to interchange the $\inf$ and $\ex$ operations.
The estimate \eqref{eq:gap-ergodic} then follows immediately.

To establish the almost sure bound \eqref{eq:gap-ergodic-as}, set $S_{\run+1} = \sum_{\iRun=\start}^{\run} \step_{\iRun} \snoise_{\iRun+1}$ and $R_{\run+1} = (2\strong)^{-1} \sum_{\iRun=\start}^{\run} \step_{\iRun}^{2} \dnorm{\est_{\iRun+1}}^{2}$.
Then, \eqref{eq:gap-ergodic2} becomes
\begin{equation}
\label{eq:gap-ergodic5}
\sum_{\iRun=\start}^{\run} \step_{\iRun} \braket{\payv(\act_{\iRun})}{\eq - \act_{\iRun}}
	\leq \fench(\eq,\score_{\start}) + S_{\run} + R_{\run},
\end{equation}
Arguing as in the proof of \cref{thm:local-imperfect}, it follows that $\sup_{\run} \exof{\abs{S_{\run}}}$ and $\sup_{\run} \exof{R_{\run}}$ are both finite, i.e. $S_{\run}$ and $R_{\run}$ are both bounded in $L^{1}$.
By Doob's (sub)martingale convergence theorem \citep[Theorem 2.5]{HH80}, it also follows that $S_{\run}$ and $R_{\run}$ both converge to an \as finite limit $S_{\infty}$ and $R_{\infty}$ respectively.
Consequently, by \eqref{eq:gap-ergodic5}, there exists a finite \as random variable $A>0$ such that
\begin{equation}
\label{eq:gap-ergodic6}
\sum_{\iRun=\start}^{\run} \step_{\iRun} \braket{\payv(\act_{\iRun})}{\eq - \act_{\iRun}}
	\leq A
	\quad
	\text{for all $\run$ \as}.
\end{equation}
The bound \eqref{eq:gap-ergodic-as} follows by taking the minimum of \eqref{eq:gap-ergodic6} over $\eq\in\eqset$ and dividing both sides by $\sum_{\iRun=\start}^{\run} \step_{\iRun}$.
Finally, applying Doob's maximal inequality to \eqref{eq:prob-bound1} and \eqref{eq:prob-bound2}, we obtain $\probof[\big]{\sup_{\run} S_{\run} \geq \noisedev \norm{\feas}\, \kappa} \leq \delta/2$ and $\probof[\big]{\sup_{\run} R_{\run} \geq \vbound^{2} \kappa^{2}} \leq \delta/2$.
Combining these bounds with \eqref{eq:gap-ergodic5} shows that $A$ can be taken to satisfy \eqref{eq:gap-constant} with probability at least $1-\delta$, as claimed.
\end{Proof}

\begin{Proof}[Proof of \cref{cor:gap-ergodic-opt}]
By the definition \eqref{eq:Fench-set} of the setwise Fenchel coupling, we have $\fench_{\start} \leq h(\eq) + h^{\ast}(0) \leq \max h - \min h = \depth$.
Our claim then follows
by invoking Jensen's inequality,
noting that $\exof{\dist(\eqset,\act_{\run})}^{2} \leq \exof{\dist(\eqset,\act_{\run})^{2}} \leq L^{-1} \exof{\gap(\act_{\run})}$,
and applying \eqref{eq:gap-ergodic}.
\end{Proof}

Although the mean bound \eqref{eq:gap-ergodic} is valid for any step-size sequence, the summability condition $\sum_{\run=\start}^{\infty} \step_{\run}^{2} < \infty$ for the almost sure bound \eqref{eq:gap-ergodic-as} rules out more aggressive step-size policies of the form $\step_{\run} \propto 1/n^{\beta}$ for $\beta \leq 1/2$.
Specifically, the ``critical'' value $\beta = 1/2$ is again tied to the finite \acl{MSE} hypothesis \eqref{eq:MSE}:
if the players' gradient measurements have finite moments up to some order $q>2$, a more refined application of Doob's inequality reveals that \eqref{eq:gap-ergodic-as} still holds under the lighter summability requirement $\sum_{\run=\start}^{\infty} \step_{\run}^{1+q/2} < \infty$.
In this case, the exponent $\beta=1/2$ is optimal with respect to the guarantee \eqref{eq:gap-ergodic} and leads to an almost sure convergence rate of the order of $\bigoh(\run^{-1/2}\log\run)$.

Except for this $\log n$ factor, the $\bigoh(\run^{-1/2})$ convergence rate of \eqref{eq:DA} is the exact lower complexity bound for black-box subgradient schemes for convex problems \citep{NY83,Nes04}.
Thus, running \eqref{eq:DA} with a step-size policy of the form $\step_{\run} \propto \run^{-1/2}$ leads to a convergence speed that is optimal in the mean, and near-optimal with high probability.
It is also worth noting that, when the horizon of play is known in advance (as in \cref{cor:gap-ergodic-opt}), the constant $\depth = \max h -\min h$ that results from the initialization $\score_{\start} = 0$ is essentially the same as the constant that appears in the stochastic \acl{MD} analysis of \cite{NJLS09} and \cite{Nes09}.

\subsection{Running length}
\label{sec:length}

Intuitively, the main obstacle to achieving rapid convergence is that, even with an optimized step-size policy, the sequence of play may end up oscillating around an equilibrium state because of the noise in the players' observations.
To study such phenomena, we focus below on the \emph{running length} of \eqref{eq:DA}, defined as
\begin{equation}
\label{eq:length}
\length_{\run}
	= \sum_{\iRun=\start}^{\run-1} \norm{\act_{\iRun+1} - \act_{\iRun}}.
\end{equation}
Obviously, if $\act_{\run}$ converges to some $\eq\in\feas$, a shorter length signifies less oscillations of $\act_{\run}$ around $\eq$.
Thus, in a certain way, $\length_{\run}$ is a more refined convergence criterion than the induced equilibrium gap $\gap(\act_{\run})$.

Our next result shows that the mean running length of \eqref{eq:DA} until players reach an $\eps$-neighborhood of a (strongly) stable set is at most $\bigoh(1/\eps^{2})$:

\begin{theorem}
\label{thm:length}
Suppose that \eqref{eq:DA} is run with
imperfect feedback satisfying \crefrange{eq:zeromean}{eq:MSE}
and
a step-size $\step_{\run}$ such that
$\sum_{\run=\start}^{\infty} \step_{\run}^{2} < \infty$ and $\sum_{\run=\start}^{\infty} \step_{\run} = \infty$.
Also, given a closed subset $\eqset$ of $\feas$, consider the stopping time $\run_{\eps} = \inf\setdef{\run\geq0}{\dist(\eqset,\act_{\run}) \leq \eps}$ and let  $\length_{\eps} \equiv \length_{\run_{\eps}}$ denote the running length of \eqref{eq:DA} until $\act_{\run}$ reaches an $\eps$-neighborhood of $\eqset$.
If $\eqset$ is $L$-strongly stable, we have
\begin{equation}
\label{eq:length-bound}
\exof{\length_{\eps}}
	\leq \frac{\vbound}{\strong L} \frac{\fench_{\start} + (2\strong)^{-1} \vbound^{2} \sum_{\iRun=\start}^{\infty} \step_{\iRun}^{2}}{\eps^{2}}.
\end{equation}
\end{theorem}

\begin{Proof}
For all $\eq\in\eqset$ and all $\run\in\N$, \eqref{eq:Dbound-stoch} yields
\begin{flalign}
\label{eq:length-bound1}
\fench(\eq,\score_{\run_{\eps}\wedge \run+1})
	&\leq \fench(\eq,\score_{\start})
	- \sum_{\iRun=\start}^{\run_{\eps}\wedge \run} \step_{\iRun} \braket{\payv(\act_{\iRun})}{\act_{\iRun} - \eq}
	\notag\\
	&+ \sum_{\iRun=\start}^{\run_{\eps}\wedge \run} \step_{\iRun} \snoise_{\iRun+1}
	+ \frac{1}{2\strong} \sum_{\iRun=\start}^{\run_{\eps}\wedge \run} \step_{\iRun}^{2} \dnorm{\est_{\iRun+1}}^{2}.
\end{flalign}
Hence, after taking expectations and minimizing over $\eq\in\eqset$, we get
\begin{equation}
\label{eq:length-bound2}
0
	\leq \fench_{\start}
	- L\eps^{2} \exof*{\sum_{\iRun=\start}^{\run_{\eps}\wedge \run} \step_{\iRun}}
	+ \exof*{\sum_{\iRun=\start}^{\run_{\eps}\wedge \run} \step_{\iRun} \snoise_{\iRun+1}}
	+ \frac{\vbound^{2}}{2\strong} \sum_{\iRun=\start}^{\infty} \step_{\iRun}^{2},
\end{equation}
where we we used the fact that $\norm{\act_{\iRun} - \eq} \geq \eps$ for all $\iRun\leq \run_{\eps}$.

Consider now
the stopped process $S_{\run_{\eps}\wedge \run} = \sum_{\iRun=\start}^{\run_{\eps}\wedge \run} \step_{\iRun} \snoise_{\iRun+1}$.
Since $\run_{\eps}\wedge \run \leq \run < \infty$, $S_{\run_{\eps}\wedge \run}$ is a martingale and $\exof{S_{\run_{\eps}\wedge \run}} = 0$.
Thus, by rearranging \eqref{eq:length-bound2}, we obtain
\begin{equation}
\label{eq:length-bound3}
\exof*{\sum_{\iRun=\start}^{\run_{\eps}\wedge \run} \step_{\iRun}}
	\leq \frac{\fench_{\start} + (2\strong)^{-1} \vbound^{2} \sum_{\iRun=\start}^{\infty} \step_{\iRun}^{2}}{L\eps^{2}}.
\end{equation}
Hence, with $\run_{\eps}\wedge\run \to \run_{\eps}$ as $\run\to\infty$, Lebesgue's monotone convergence theorem shows that the process $\tau_{\eps} = \sum_{\iRun=\start}^{\run_{\eps}} \step_{\iRun}$ is finite in expectation and
\begin{equation}
\label{eq:length-bound4}
\exof{\tau_{\eps}}
	\leq \frac{\fench_{\start} + (2\strong)^{-1} \vbound^{2} \sum_{\iRun=\start}^{\infty} \step_{\iRun}^{2}}{L \eps^{2}}.
\end{equation}

Furthermore, by \cref{prop:choice} and the definition of $\length_{\run}$, we also have
\begin{equation}
\label{eq:length-bound5}
\length_{\run}
	= \sum_{\iRun=\start}^{\run-1} \norm{\act_{\iRun+1} - \act_{\iRun}}
	\leq \frac{1}{\strong} \sum_{\iRun=\start}^{\run-1} \dnorm{\score_{\iRun} - \score_{\iRun-1}}
	= \frac{1}{\strong} \sum_{\iRun=\start}^{\run-1} \step_{\iRun} \dnorm{\est_{\iRun+1}}.
\end{equation}
Now,
let $\zeta_{\iRun+1} = \dnorm{\est_{\iRun+1}}$
and $\Psi_{\run+1} = \sum_{\iRun=\start}^{\run} \step_{\iRun} \left[ \zeta_{\iRun+1} - \exof{\zeta_{\iRun+1} \given \filter_{\iRun}}\right]$.
By construction, $\Psi_{\run}$ is a martingale and
\begin{equation}
\exof{\Psi_{\run+1}^{2}}
	= \exof*{\sum_{\iRun=\start}^{\run} \step_{\iRun}^{2} \left[ \zeta_{\iRun+1} - \exof{\zeta_{\iRun+1} \given \filter_{\iRun}}\right]^{2}} \leq 2 \vbound^{2} \sum_{\iRun=\start}^{\infty} \step_{\iRun}^{2} < \infty
	\quad
	\text{for all $\run$}.
\end{equation}
Thus, by the optional stopping theorem \citep[p.~485]{Shi95}, we get $\exof{\Psi_{\run_{\eps}}} = \exof{\Psi_{\start}} = 0$, so
\begin{equation}
\label{eq:length-stopping}
\exof*{\sum_{\iRun=\start}^{\run_{\eps}} \step_{\iRun} \zeta_{\iRun+1}}
	= \exof*{\sum_{\iRun=\start}^{\run_{\eps}} \step_{\iRun} \exof{\zeta_{\iRun+1} \given \filter_{\iRun}}}
	\leq \vbound \exof*{\sum_{\iRun=\start}^{\run_{\eps}} \step_{\iRun}}
	= \vbound \exof{\tau_{\eps}}.
\end{equation}
Our claim then follows by combining \eqref{eq:length-bound5} and \eqref{eq:length-stopping} with the bound \eqref{eq:length-bound4}.
\end{Proof}

\cref{thm:length} should be contrasted to classic results on the Kurdyka\textendash \L ojasiewicz inequality where having a ``bounded length'' property is crucial in establishing trajectory convergence \citep{BDLM10}.
In our stochastic setting, it is not realistic to expect a bounded length (even on average), because, generically, the noise does not vanish in the neighborhood of a \acl{NE}.%
\footnote{For a notable exception however, see \cref{thm:rate-sharp} below.}
Instead, \cref{thm:length} should be interpreted as a measure of how the fluctuations due to noise and uncertainty affect the trajectories' average length;
the authors are not aware of any similar results along these lines.

\subsection{Sharp equilibria and fast convergence}
\label{sec:rate-sharp}

Because of the random shocks induced by the noise in the players' gradient observations, it is difficult to obtain an almost sure (or high probability) estimate for the convergence rate of the last iterate $\act_{\run}$ of \eqref{eq:DA}.
Specifically,
even with a rapidly decreasing step-size policy,
a single realization of the error process $\noise_{\run}$ may lead to an arbitrarily big jump of $\act_{\run}$ at any time,
thus destroying any almost sure bound on the convergence rate of $\act_{\run}$.

On the other hand, in finite games, \cite{CHM17-SAGT} recently showed that logit-based learning (cf. \cref{alg:logit}) achieves a quasi-linear convergence rate with high probability if the equilibrium in question is strict.
Specifically, \cite{CHM17-SAGT} showed that if $\eq$ is a strict \acl{NE} and $\act_{\run}$ does not start too far from $\eq$, then, with high probability, $\norm{\act_{\run} - \eq} = \bigoh(-\sharp \sum_{\iRun=\start}^{\run} \step_{\iRun})$ for some positive constant $\sharp>0$ that depends only on the players' relative payoff differences.

Building on the variational characterization of strict \aclp{NE} provided by \cref{prop:strict}, we consider below the following analogue for continuous games:

\begin{definition}
\label{def:sharp}
We say that $\eq\in\feas$ is a \emph{sharp equilibrium} of $\game$ if
\begin{equation}
\label{eq:sharp}
\braket{\payv(\eq)}{z}
	\leq 0
	\quad
	\text{for all $z\in\tcone(\eq)$},
\end{equation}
with equality if and only if $z=0$.
\end{definition}

\begin{remark}
The terminology ``sharp'' follows \citet[Chapter~5.2]{Pol87}, who introduced a similar notion for (unconstrained) convex programs.
In particular, in the single-player case, it is easy to see that \eqref{eq:sharp} implies that $\eq$ is a \emph{sharp maximum} of $\pay(x)$, i.e. $\pay(\eq) - \pay(x) \geq \sharp\norm{x-\eq}$ for some $\sharp>0$.
\end{remark}

A first consequence of \cref{def:sharp} is that $\payv(\eq)$ lies in the topological interior of the polar cone $\pcone(\eq)$ to $\feas$ at $\eq$ (for a schematic illustration, see \cref{fig:Nash});
in turn, this implies that sharp equilibria can only occur at \emph{corners} of $\feas$.
By continuity, this further implies that sharp equilibria are locally stable (cf. the proof of \cref{thm:rate-sharp} below);
hence, by \cref{prop:Nash-stable}, sharp equilibria are also isolated.
Our next result shows that if players employ \eqref{eq:DA} with surjective choice maps, then, with high probability, sharp equilibria are attained in a \emph{finite} number of steps:

\begin{theorem}
\label{thm:rate-sharp}
Fix a tolerance level $\delta>0$ and suppose that \eqref{eq:DA} is run with surjective choice maps and a sufficiently small step-size $\step_{\run}$ such that
$\sum_{\run=\start}^{\infty} \step_{\run}^{2} < \infty$ and $\sum_{\run=\start}^{\infty} \step_{\run} = \infty$.
If $\eq$ is sharp and \eqref{eq:DA} is not initialized too far from $\eq$, we have
\begin{equation}
\label{eq:rate-sharp}
\probof{\textup{$\act_{\run}$ reaches $\eq$ in a \emph{finite} number of steps}}
	\geq 1-\delta,
\end{equation}
provided that \crefrange{eq:zeromean}{eq:Lipschitz} hold.
If, in addition, $\eq$ is globally stable, $\act_{\run}$ converges to $\eq$ in a finite number of steps from every initial condition \as.
\end{theorem}

\begin{Proof}
As we noted above, $\payv(\eq)$ lies in the interior of the polar cone $\pcone(\eq)$ to $\feas$ at $\eq$.%
\footnote{Indeed, if this were not the case, we would have $\braket{\payv(\eq)}{z} = 0$ for some nonzero $z\in\tcone(\eq)$.}
Hence, by continuity, there exists a neighborhood $\eqnhd$ of $\eq$ such that $\payv(x)\in\intr(\pcone(\eq))$ for all $x\in \eqnhd$.
In turn, this implies that $\braket{\payv(x)}{x - \eq} < 0$ for all $x\in \eqnhd \setminus\{\eq\}$, i.e. $\eq$ is stable.
Therefore, by \cref{thm:local-imperfect}, there exists a neighborhood $U$ of $\eq$ such that $\act_{\run}$ converges to $\eq$ with probability at least $1-\delta$.

Now, let $U'\subseteq \eqnhd$ be a sufficiently small neighborhood of $\eq$ such that $\braket{\payv(x)}{z} \leq - \sharp\norm{z}$ for some $\sharp>0$ and for all $z\in\tcone(\eq)$.%
\footnote{That such a neighborhood exists is a direct consequence of \cref{def:sharp}.}
Then, with probability at least $1-\delta$, there exists some (random) $\run_{0}$ such that $\act_{\run} \in U'$ for all $\run\geq\run_{0}$, so $\braket{\payv(\act_{\run})}{z} \leq -\sharp\norm{z}$ for all $\run\geq \run_{0}$.
Thus, for all $z\in\tcone(\eq)$ with $\norm{z} = 1$, we have
\begin{flalign}
\label{eq:sharp-bound1}
\braket{\score_{\run+1}}{z}
	&= \braket{\score_{\run_{0}}}{z}
	+ \sum_{\iRun=\run_{0}}^{\run} \step_{\iRun} \braket{\payv(\act_{\iRun})}{z}
	+ \sum_{\iRun=\run_{0}}^{\run} \step_{\iRun} \braket{\noise_{\iRun+1}}{z}
	\notag\\
	&\leq \dnorm{\score_{\run_{0}}}
	- \sharp \sum_{\iRun=\run_{0}}^{\run} \step_{\iRun}
	+ \sum_{\iRun=\run_{0}}^{\run} \step_{\iRun} \braket{\noise_{\iRun+1}}{z}.
\end{flalign}
By the \acl{LLN} for \aclp{MDS} \citep[Theorem 2.18]{HH80}, we also have $\sum_{\iRun=\run_{0}}^{\run} \step_{\iRun} \noise_{\iRun+1} / \sum_{\iRun=\run_{0}}^{\run} \step_{\iRun} \to 0$ \as, so there exists some $\run^{\ast}$ such that $\dnorm{\sum_{\iRun=\run_{0}}^{\run} \step_{\iRun} \noise_{\iRun+1}} \leq (\sharp/2) \sum_{\iRun=\run_{0}}^{\run} \step_{\iRun}$ for all $\run\geq\run^{\ast}$ \as.
We thus obtain
\begin{equation}
\label{eq:sharp-bound2}
\braket{\score_{\run+1}}{z}
	\leq \dnorm{\score_{\run_{0}}}
	- \sharp \sum_{\iRun=\run_{0}}^{\run} \step_{\iRun}
	+ \frac{\sharp}{2} \norm{z} \sum_{\iRun=\run_{0}}^{\run} \step_{\iRun}
	\leq \dnorm{\score_{\run_{0}}}
	- \frac{\sharp}{2} \sum_{\iRun=\run_{0}}^{\run} \step_{\iRun},
\end{equation}
showing that $\braket{\score_{\run}}{z}\to -\infty$ uniformly in $z$ with probability at least $1-\delta$.

To proceed, \cref{prop:preimage} in \cref{app:aux} shows that $\deq + \pcone(\eq) \subseteq \choice^{-1}(\eq)$ whenever 
$\choice(\deq) = \eq$.
Since $\choice$ is surjective, there exists some $\deq\in\choice^{-1}(\eq)$, so it suffices to show that, with probability at least $1-\delta$, $\score_{\run}$ lies in the pointed cone $\deq + \pcone(\eq)$ for all sufficiently large $\run$.
To do so, simply note that $\score_{\run} - \deq \in \pcone(\eq)$ if and only if $\braket{\score_{\run} - \deq}{z} \leq 0$ for all $z\in\tcone(\eq)$ with $\norm{z}=1$.
Since $\braket{\score_{\run}}{z}$ converges uniformly to $-\infty$ with probability at least $1-\delta$, our assertion is immediate.

Finally, for the globally stable case, recall that $\act_{\run}$ converges to $\eq$ with probability $1$ from any initial condition (\cref{thm:global-imperfect}).
The argument above shows that $\act_{\run} = \eq$ for all large $\run$, so $\act_{\run}$ converges to $\eq$ in a finite number of steps \as.
\end{Proof}

\begin{remark}
\cref{thm:rate-sharp} suggests that \acl{DA} with surjective choice maps leads to significantly faster convergence to sharp equilibria.
In this way, it is consistent with an observation made by \citet[Proposition 5.2]{MS16} for the convergence of the continuous-time, deterministic dynamics \eqref{eq:DA-cont} in finite games.
\end{remark}

%% file: Discussion.tex

An important question in the implementation of \acl{DA} is the choice of regularizer, which in turn determines the players' choice maps $\choice_{\play}\from\dual_{\play}\to\feas_{\play}$.
From a qualitative point of view, this choice would not seem to matter much:
the convergence results of \cref{sec:analysis,sec:finite} hold for all choice maps of the form \eqref{eq:choice}.
Quantitatively however, the specific choice map employed by each player impacts the algorithm's convergence speed, and different choice maps could lead to vastly different rates of convergence.

As noted above, in the case of sharp equilibria, this choice seems to favor nonsteep penalty functions (that is, surjective choice maps).
Nonetheless, in the general case, the situation is less clear because of the dimensional dependence hidden in the $\depth/\strong$ factor that appears e.g. in the mean rate guarantee \eqref{eq:gap-ergodic-opt}.
This factor depends crucially on the geometry of the players' action spaces and the underlying norm, and its optimum value may be attained by \emph{steep} penalty functions \textendash\ for instance, the entropic regularizer \eqref{eq:entropy} is well known to be asymptotically optimal in the case of simplex-like feasible regions \citefp[p.~140]{SS11}.

Another key question in game-theoretic and online learning has to do with the information that is available to the players at each stage.
If players perform a two-point sampling step in order to simulate an extra oracle call at an action profile different than the one employed, this extra information could be presumably leveraged in order to increase the speed of convergence to a \acl{NE}.
In an offline setting, this can be achieved by more sophisticated techniques relying on dual extrapolation \citefp{Nes07} and/or mirror-prox methods \citefp{JNT11}.
Extending these extra-gradient approaches to online learning processes as above would be an interesting extension of the current work.

At the other end of the spectrum, if players only have access to their realized, in-game payoffs, they would need to reconstruct their individual payoff gradients via a suitable single-shot estimator \citep{Pol87,FKM05}.
We believe our convergence analysis can be extended to this case by properly controlling the ``bias-variance'' tradeoff of this estimator and using more refined stochastic approximation arguments.
The very recent manuscript by \citef{BBF16} provides an encouraging first step in the case of (strictly) concave games with one-dimensional action sets;
we intend to explore this direction in future work.

%% file: App-Auxiliary.tex

\newcommand{\region}{\mathcal{R}}

In this appendix, we collect some auxiliary results that would have otherwise disrupted the flow of the main text.
We begin with the basic properties of the Fenchel coupling:

\begin{Proof}[Proof of \cref{prop:Fenchel}]
For our first claim, let $x=\choice(y)$.
Then, by definition
\begin{equation}
\label{eq:Fench-sub}
\fench(\base,y)
	= h(\base) + \braket{y}{\choice(y)} - h(\choice(y)) - \braket{y}{\base}
	= h(\base) - h(x) - \braket{y}{\base - x}.
\end{equation}
Since $y\in\pd h(x)$ by \cref{prop:choice}, we have $\braket{y}{\base - x} = h'(x;\base-x)$ whenever $x\in\intacts$, thus proving \eqref{eq:Fench-Bregman}.
Furthermore, the strong convexity of $h$ also yields
\begin{flalign}
h(x) + t \braket{y}{\base - x}
	&\leq h(x + t(\base - x))
	\notag\\
	&\leq th(\base) + (1-t)h(x) - \tfrac{1}{2} \strong t(1-t) \norm{x - \base}^{2},
\end{flalign}
leading to the bound
\begin{equation}
\label{eq:divbound}
\tfrac{1}{2} \strong(1-t) \norm{x - \base}^{2}
	\leq h(\base) - h(x) - \braket{y}{\base - x}
	= \fench(\base,y)
\end{equation}
for all $t\in(0,1]$.
\cref{eq:Fench-norm} then follows by letting $t\to0^{+}$ in \eqref{eq:divbound}.

Finally, for our third claim, we have
\begin{flalign}
\fench(\base,y')
	&= h(\base) + h^{\ast}(y') - \braket{y'}{\base}
	\notag\\
	&\leq h(\base) + h^{\ast}(y) + \braket{y' - y}{\nabla h^{\ast}(y)} + \frac{1}{2\strong} \dnorm{y' - y}^{2} - \braket{y'}{\base}
	\notag\\
	&= \fench(\base,y) + \braket{y' - y}{\choice(y) - \base} + \frac{1}{2\strong} \dnorm{y' - y}^{2},
\end{flalign}
where the inequality in the second line follows from the fact that $h^{\ast}$ is $(1/\strong)$-strongly smooth \citefp[Theorem 12.60(e)]{RW98}.
\end{Proof}

Complementing \cref{prop:Fenchel}, our next result concerns the inverse images of the choice map $\choice$:

\begin{proposition}
\label{prop:preimage}
Let $h$ be a penalty function on $\feas$,
and let $\eq\in\feas$.
If $\eq = \choice(\deq)$ for some $\deq\in\dual$, then $\deq + \pcone(\eq) \subseteq \choice^{-1}(\eq)$.
\end{proposition}

\begin{Proof}
By \cref{prop:choice}, we have $\eq = \choice(y)$ if and only if $y\in\pd h(\eq)$, so it suffices to show that $\deq + \payv \in \pd h(\eq)$ for all $\payv\in\pcone(\eq)$.
Indeed, we have $\braket{\payv}{x-\eq} \leq 0$ for all $x\in\feas$, so
\begin{equation}
h(x)
	\geq h(\eq) + \braket{\deq}{x - \eq}
	\geq h(\eq) + \braket{\deq + \payv}{x - \eq}.
\end{equation}
The above shows that $\deq + \payv \in \pd h(\eq)$, as claimed.
\end{Proof}

Our next result concerns the evolution of the Fenchel coupling under the dynamics \eqref{eq:DA-cont}:

\begin{lemma}
\label{lem:Lyapunov}
Let $x(t) = \choice(y(t))$ be a solution orbit of \eqref{eq:DA-cont}.
Then, for all $\base\in\feas$, we have
\begin{equation}
\label{eq:Lyapunov}
\frac{d}{dt} \fench(\base,y(t))
	= \braket{\payv(x(t))}{x(t) - \base}.
\end{equation}
\end{lemma}

\begin{Proof}
By definition, we have
\begin{flalign}
\label{eq:dF}
\frac{d}{dt} \fench(\base,y(t))
	&= \frac{d}{dt} \left[ h(\base) + h^{\ast}(y(t)) - \braket{y(t)}{\base} \right]
	\notag\\
	&= \braket{\dot y(t)}{\nabla h^{\ast}(y(t))} - \braket{\dot y(t)}{\base}
	= \braket{\payv(x(t))}{x(t) - \base},
\end{flalign}
where, in the last line, we used \cref{prop:choice}.
\end{Proof}

Our last auxiliary result shows that, if the sequence of play generated by \eqref{eq:DA} is contained in the ``basin of attraction'' of a stable set $\eqset$, then it admits an accumulation point in $\eqset$:

\begin{lemma}
\label{lem:visit}
Suppose that $\eqset\subseteq\feas$ is stable and \eqref{eq:DA} is run with a step-size such that
$\sum_{\run=\start}^{\infty} \step_{\run}^{2} < \infty$ and $\sum_{\run=\start}^{\infty} \step_{\run} = \infty$.
Assume further that $(\act_{\run})_{\run=\start}^{\infty}$ is contained in a region $\region$ of $\feas$ such that \eqref{eq:VS} holds for all $x\in\region$.
Then, under \labelcref{eq:zeromean,eq:MSE}, every neighborhood $U$ of $\eqset$ is recurrent;
specifically, there exists a subsequence $\act_{\run_{\iRun}}$ of $\act_{\run}$ such that $\act_{\run_{\iRun}}\to\eqset$ \as.
Finally, if \eqref{eq:DA} is run with perfect feedback \textpar{$\noisedev=0$}, the above holds under the lighter assumption $\sum_{\iRun=\start}^{\run} \step_{\iRun}^{2} \big/ \sum_{\iRun=\start}^{\run} \step_{\iRun} \to 0$.
\end{lemma}

\begin{Proof}[Proof of \cref{lem:visit}]
Let $U$ be a neighborhood of $\eqset$ and assume to the contrary that, with positive probability, $\act_{\run}\notin U$ for all sufficiently large $\run$.
By starting the sequence at a later index if necessary, we may assume that $\act_{\run}\notin U$ for all $\run$ without loss of generality.
Thus, with $\eqset$ stable and $\act_{\run}\in\region$ for all $\run$ by assumption, there exists some $\sharp>0$ such that
\begin{equation}
\label{eq:drift-bound}
\braket{\payv(\act_{\run})}{\act_{\run} - \eq}
	\leq -\sharp
	\quad
	\text{for all $\eq\in\eqset$ and for all $\run$}.
\end{equation}
As a result, for all $\eq\in\eqset$, we get
\begin{flalign}
\label{eq:Dbound-1}
\fench(\eq,\score_{\run+1})
	&= \fench(\eq,\score_{\run} + \step_{\run} \est_{\run+1})
	\notag\\
	&\leq \fench(\eq,\score_{\run})
	+ \step_{\run} \braket{\payv(\act_{\run}) + \noise_{\run+1}}{\act_{\run} - \eq}
	+ \frac{1}{2\strong} \step_{\run}^{2} \dnorm{\est_{\run+1}}^{2}
	\notag\\
	&\leq \fench(\eq,\score_{\run})
	- \sharp\step_{\run}
	+ \step_{\run} \snoise_{\run+1}
	+ \frac{1}{2\strong} \step_{\run}^{2} \dnorm{\est_{\run+1}}^{2},
\end{flalign}
where we used \cref{prop:Fenchel} in the second line and we set $\snoise_{\run+1} = \braket{\noise_{\run+1}}{\act_{\run} - \eq}$ in the third.
Telescoping \eqref{eq:Dbound-1} then gives
\begin{equation}
\label{eq:Dbound-2}
\fench(\eq,\score_{\run+1})
	\leq \fench(\eq,\score_{\start})
	- \tau_{\run} \left[
	\sharp
	- \frac{\sum_{\iRun=\start}^{\run} \step_{\iRun} \snoise_{\iRun+1}}{\tau_{\run}}
	- \frac{1}{2\strong} \frac{\sum_{\iRun=\start}^{\run} \step_{\iRun}^{2} \dnorm{\est_{\iRun+1}}^{2}}{\tau_{\run}}
	\right],
\end{equation}
where $\tau_{\run} = \sum_{\iRun=\start}^{\run} \step_{\iRun}$.

Since $\exof{\snoise_{\run+1} \given \filter_{\run}} = \braket{\exof{\noise_{\run+1} \given \filter_{\run}}}{\act_{\run} - \eq} = 0$ by \eqref{eq:zeromean} and $\exof{\abs{\snoise_{\run+1}}^{2} \given \filter_{\run}} \leq \exof{\dnorm{\noise_{\run+1}}^{2} \norm{\act_{\run} - \eq}^{2} \given \filter_{\run}} \leq \noisevar \, \norm{\feas}^{2} < \infty$ by \eqref{eq:MSE}, the \acl{LLN} for \aclp{MDS} yields $\tau_{\run}^{-1} \sum_{\iRun=\start}^{\run} \step_{\iRun} \snoise_{\iRun+1} \to 0$ \citefp[Theorem 2.18]{HH80}.
Furthermore, letting $R_{\run+1} = \sum_{\iRun=\start}^{\run} \step_{\iRun}^{2} \dnorm{\est_{\iRun+1}}^{2}$, we also get
\begin{equation}
\exof{R_{\run+1}}
	\leq \sum_{\iRun=\start}^{\run} \step_{\iRun}^{2} \exof{\est_{\iRun+1}}^{2}
	\leq \vbound^{2} \sum_{\iRun=\start}^{\infty} \step_{\iRun}^{2}
	< \infty
	\quad
	\text{for all $\run$},
\end{equation}
so Doob's martingale convergence theorem shows that $R_{\run}$ converges \as to some random, finite value \citefp[Theorem 2.5]{HH80}.

Combining the above, \eqref{eq:Dbound-2} gives $\fench(\eq,\score_{\run}) \sim -a \tau_{\run} \to -\infty$ \as, a contradiction.
Finally, if $\noisedev=0$, we also have $\snoise_{\run+1} = 0$ and $\dnorm{\est_{\run+1}}^{2} = \dnorm{\payv(\act_{\run})}^{2} \leq \vbound^{2}$ for all $\run$, so \eqref{eq:Dbound-2} yields $\fench(\eq,\score_{\run})\to-\infty$ provided that $\tau_{\run}^{-1} \sum_{\iRun=\start}^{\run} \step_{\iRun}^{2} \to 0$, a contradiction.
In both cases, we conclude that $\act_{\run}$ is recurrent, as claimed.
\end{Proof}

Finally, we turn to the characterization of strict equilibria in finite games:

\begin{Proof}[Proof of \cref{prop:strict}]
We will show that (\emph{a})$\implies$(\emph{b})$\implies$(\emph{c})$\implies$(\emph{a}).

\paragraph{\upshape(\emph{a})$\implies$(\emph{b})}
Suppose that $\eq = (\peq_{1},\dotsc,\peq_{\nPlayers})$ is a strict equilibrium.
Then, the weak inequality $\braket{\payv(\eq)}{z} \leq 0$ follows from \cref{prop:Nash-var}.
For the strict part, if $z_{\play}\in\tcone_{\play}(\eq_{\play})$ is nonzero for some $\play\in\players$, we readily get
\begin{equation}
\braket{\payv_{\play}(\eq)}{z_{\play}}
	= \sum_{\pure_{\play}\neq\peq_{\play}} z_{\play,\pure_{\play}}
	\left[
	\pay_{\play}(\peq_{\play};\peq_{-\play}) - \pay_{\play}(\pure_{\play};\peq_{-\play})
	\right]
	< 0,
\end{equation}
where we used the fact that $z_{\play}$ is tangent to $\feas$ at $\eq_{\play}$, so $\sum_{\pure_{\play}\in\pures_{\play}} z_{\play\pure_{\play}} = 0$ and $z_{\play\pure_{\play}} \geq 0$ for $\pure_{\play}\neq\peq_{\play}$, with at least one of these inequalities being strict when $z_{\play}\neq0$.

\paragraph{\upshape(\emph{b})$\implies$(\emph{c})}
Property (\emph{b}) implies that $\payv(\eq)$ lies in the interior of the polar cone $\pcone(\eq)$ to $\feas$ at $\eq$.
Since $\pcone(\eq)$ has nonempty interior, continuity implies that $\payv(x)$ also lies in $\pcone(\eq)$ for $x$ sufficiently close to $\eq$.
We thus get $\braket{\payv(x)}{x - \eq} \leq 0$ for all $x$ in a neighborhood of $\eq$, i.e. $\eq$ is stable.

\smallskip

\paragraph{\upshape(\emph{c})$\implies$(\emph{a})}
Assume that $\eq$ is stable but not strict, so $\pay_{\play\pure_{\play}}(\eq) = \pay_{\play\purealt_{\play}}(\eq)$ for some $\play\in\players$, and some $\pure_{\play}\in\supp(\eq_{\play})$, $\purealt_{\play}\in\pures_{\play}$.
Then, if we take $x_{\play} = \eq_{\play} + \lambda (\bvec_{\play\purealt_{\play}} - \bvec_{\play\pure_{\play}})$ and $x_{-\play} = \eq_{-\play}$ with $\lambda>0$ small enough, we get
\begin{equation}
\braket{\payv(x)}{x - \eq}
	= \braket{\payv_{\play}(x)}{x_{\play} - \eq_{\play}}
	= \lambda \pay_{\play\purealt_{\play}}(\eq) - \lambda \pay_{\play\pure_{\play}}(\eq)
	= 0,
\end{equation}
contradicting the assumption that $\eq$ is stable.
This shows that $\eq$ is strict.
\end{Proof}